\documentclass{article}

\usepackage[all]{nowidow}
\usepackage[protrusion=true, final]{microtype}
\usepackage{josealvim}
\usepackage{hyperref}
\usepackage[backend=biber]{biblatex}
\usepackage{quiver}
\usepackage[inline]{enumitem}
\usepackage{import}
\usepackage{verbatim}
\usepackage[utf8]{inputenc}

\begin{document}
	\title{
		-- \(\QSet\) \textit{\&} Friends --
		\\
		Regarding Singleton and Gluing Completeness
	}
	\author{
        José G. Alvim\email{alvim@ime.usp.br}
		\and 
        Hugo L. Mariano\email{hugomar@ime.usp.br}
		\and
        Caio de A. Mendes\email{caio.mendes@usp.br}
	}
	\begin{titlepage}
	
		\maketitle
		\begin{abstract}
            This work is largely focused on extending D. Higgs' 
            \(\Omega\)-sets to the context of quantales, following
            the broad program of \cite{HOHLE199115}, we explore the
            rich category of \(\quantale Q\)-sets for strong, 
            integral and commutative quantales, or other similar 
            axioms. The focus of this work is to study the different
            notion of “completeness” a \(\quantale Q\)-set may enjoy
            and their relations, completion functors, resulting 
            reflective subcategories, their relations to relational
            morphisms. 

            We establish the general equivalence of singleton
            complete \(\quantale Q\)-sets with functional morphisms
            and the category of \(\quantale Q\)-sets with relational
            morphisms; we provide two characterizations of singleton
            completeness in categorical terms; we show that the 
            singleton complete categorical inclusion creates limits.
		\end{abstract}
		\thispagestyle{empty}
	\end{titlepage}

	\tableofcontents
    \newpage

	\section{Introduction}
\subsection*{History \textit{\&} Motivation}




The notion of sheaf on a topological space depends only on the (complete) lattice of the open sets of the space, thus it is straightforward to define sheaves for “spaces  without points”, that is, for a category \(H,\leq)\) for a given  locale \(H\) (see \cite{borceux1994handbook3}). 

In the 1970s, the topos of sheaves over a locale (= complete Heyting 
algebra) \(\mathbb{H}\) was described,
alternatively, as a category of \(\mathbb{H}\)-sets 
\cite{fourman1979sheaves}. More precisely, in 
\cite{borceux1994handbook3}, there were three categories whose objects 
were locale valued sets that are equivalent to the category of sheaves 
over a locale \(\mathbb{H}\). Two different notions of separability and 
completeness have been proposed. On the one hand, the traditional 
notions of these properties in \(\Sh(\mathbb{H})\) can be translated to 
appropriate definitions in \(\mathbb{H}\)-sets. In addition, Scott's 
notion of singletons, a definition that is inspired from the ordinary 
singleton set, leads alternative notions of completeness and 
separability. Despite some folkloric misconceptions about those 
notions, a simple counter-example (in a finite boolean algebra) shows 
that these two definitions of completeness are not logically equivalent 
-- and do not give rise to equivalent full subcategories either.

This wealth of definitions and notions, however, has had the 
unfortunate effect of muddying any discussion concerning those kinds of 
objects. A veritably deep folklore has taken root in the field which 
hinders careful thought and frightens -- almost to the point of panic 
-- anyone who is paying attention.


Nevertheless, there is a non-commutative and non-idempotent 
generalization of locales called “quantales”, introduced by C.J. Mulvey 
\cite{mulvey86quantales}. Quantales show up in logic 
\cite{yetter_1990}, and in the study of \(C^*\)-algebras 
\cite{rosenthal1990quantales}.

Many  notions of sheaf over a  quantal and  quantal-valued set are studied in many works 
(\cite{walters1981sheaves}, \cite{borceux1986quantales}, \cite{mulveyquantale}, \cite{Mulvey1995}, \cite{miraglia1998sheaves}, \cite{coniglio2000non} \cite{HOHLE199115}, \cite{heymans2012grothendieck}, 
\cite{resende2012groupoid}, \cite{zambrano}, \cite{tenorio2022sheaves}). In many cases, the base quantales  are \textit{right-sided and idempotent}.  
Herein we continue our study of quantale-valued set on 
\emph{commutative and semicartesian} quantales, initiated in \cite{AMM2023constructions}. Our approach is similar to the last one 
but, since every idempotent semicartesian quantale is a locale 
(Proposition \ref{prop:idempotent + semicartesian = locale}), our 
axioms and theirs are orthogonal in some sense.

The goal of the present work is to examine these two notions of 
completeness (and separability): 
\begin{enumerate*}[label=(\roman*)]
    \item
        via (unique) gluing of compatible families 
        (gluing completeness of \(\quantale Q\)-sets);
    \item
        via (unique) representation of (\emph{strict}) 
        singletons (Scott completeness).
\end{enumerate*}
their relations and the properties of the full subcategories of 
their examples. 

\subsection*{Main results and the paper's structure}
We have shown that: 
\begin{enumerate}
    \item 
        Scott-completeness implies gluing-completeness;
    \item   
        Both full subcategories of gluing-complete and 
        Scott-complete \(\quantale Q\)-sets with functional 
        morphisms are reflective; 
    \item 
        For “strong” \cite{HOHLE199115} quantales, it makes 
        sense to speak of “the Scott completion” of a given 
        \(\quantale Q\)-set; 
    \item
        Every \(\quantale Q\)-set is relationally isomorphic to
        its own Scott completion, completion is invariant under 
        functional isomorphisms;
    \item
        For strong quantales, the categories of Scott-complete 
        \(\quantale Q\)-sets with relational morphisms and the 
        one with morphisms are isomorphic.
    \item
        For strong quantales,
        that \(X\) being Scott-complete is equivalent to the 
        its functional representable functor “being the same” 
        as its relational representable functor.
    \item 
        For strong quantales,
        that \(X\) being Scott-complete is equivalent to its
        representable functor being “invariant” under the 
        completion endofunctor.
    \item
        For strong quantales,
        that the full inclusion of Scott-complete 
        \(\quantale Q\)-sets not only preserves (due to it 
        being a right adjoint) but also \emph{creates} limits.
        And specifically, that limits of complete 
        \(\quantale Q\)-sets can be computed “pointwise”.
\end{enumerate}





	\section{Preliminaries: Quantales}
\label{sec:quantales}

\begin{definition}
    A \emph{quantale} is a type of structure 
    \(\quantale Q = (|\quantale Q|, \leq, \tensor)\) for which 
    \((|\quantale Q|, \leq)\) is a complete lattice; 
    \((|\quantale Q|, \tensor)\) is a semigroup%
    \footnote{
        \lat{i.e.} the  binary operation 
        \(\tensor:\quantale Q\times\quantale Q\to\quantale Q\) 
        (called multiplication) is associative.
    }; and, moreover, \(\quantale Q\) is required to satisfy the 
    following distributive laws: for all \(a\in\quantale Q\) and 
    \(B\subseteq\quantale Q\),
    \begin{align*}
            a \tensor \left(\bigvee_{b\in B} b\right) 
        &=	\bigvee_{b\in B}\left(a\tensor b\right)
        \\
            \left(\bigvee_{b\in B} b\right) \tensor a 
        &=	\bigvee_{b\in B}\left(b\tensor a\right)
    \end{align*}
    We denote by \(\extent\quantale Q\) the subset of 
    \(\quantale Q\) comprised of its idempotent elements.
\end{definition}

\begin{remark}~\label{remark:quantale properties}
    \begin{enumerate}
        \item 
            In any quantale \(\quantale Q\) the 
            multiplication is increasing in both entries;
        \item 
            Since \(\bot\) is also the supremum of 
            \(\emptyset\), for any \(a\),
            \(a\tensor\bot=\bot=\bot\tensor a\)
        \item 
            Since \(\top\) is \(\sup\quantale Q\), then
            \(
                \top\tensor\top=
                \sup_{a,b}a\tensor b =
                \sup\img\tensor
            \)
    \end{enumerate}
\end{remark}

\begin{remark} 
    If \((\quantale Q, \leq)\) is a complete lattice for which the 
    binary infimum satisfies the above distributive laws, the 
    resulting quantale has \(\top\) as its unit and is -- in fact 
    -- a locale. Conversely, every locale is a unital quantale in
    such a manner.
\end{remark}

\begin{definition}
    A quantale \(\quantale Q\) is said to be
    \begin{enumerate}[label=\(\bullet\)]    
        \item \emph{bidivisible}
            when \[
                a\leq b\implies 
                \exists \lambda,\rho: a\tensor\rho = b = \lambda\tensor a
            \]
            left (right) divisibility means to drop the 
            \(\rho\) (\(\lambda\)) portion of the axiom.
        \item \emph{integral}      
            when \(\top\tensor a = a = a\tensor\top\). 
            we say it's right-sided when the right equality
            holds, and left-sided when the left equality 
            holds.
        \item \emph{unital}
            when \(\tensor\) has a unit;
        \item \emph{semicartesian} 
            when \(a \tensor b \leq a\wedge b\)
        \item \emph{commutative}
            when \(\tensor\) is;
        \item \emph{idempotent}
            when \(a \tensor a = a\);
        \item \emph{linear and strict}
            when \(\leq\) is a linear order and for \[
                a\neq \bot\ and \ a\tensor b = a\tensor c\ or\ b\tensor a = c\tensor a  \implies 
                b = c
            \]
        \item \emph{strong} 
            when for any \(e\) and \(A\)
            \cite[cf. p. 30]{HOHLE199115},
        \[
            e=e\tensor e \implies 
            e\leq\bigvee_{a\in A} a \implies 
            e\leq\bigvee_{a\in A} a \tensor a
        \]
    \end{enumerate}
    We offer the following diagram to explain some of 
    the relations between those definitions:
    
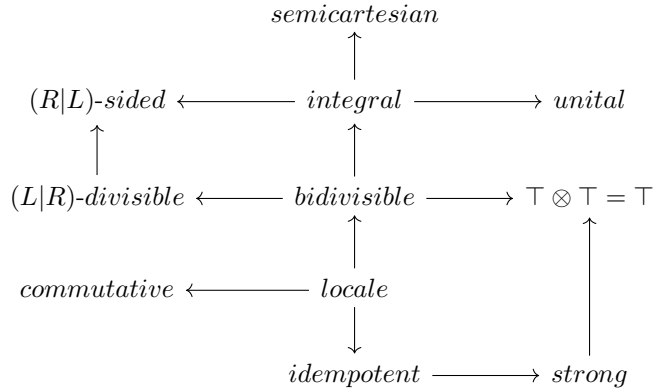
\[\begin{tikzcd}
        & semicartesian \\
        {(R|L)\text-sided} & integral & unital \\
        {(L|R)\text-divisible} & bidivisible & {\top\tensor\top=\top} \\
        commutative & locale \\
        & idempotent & strong
        \arrow[from=3-2, to=3-1]
        \arrow[from=3-2, to=2-2]
        \arrow[from=4-2, to=3-2]
        \arrow[from=4-2, to=4-1]
        \arrow[from=2-2, to=2-1]
        \arrow[from=2-2, to=1-2]
        \arrow[from=2-2, to=2-3]
        \arrow[from=4-2, to=5-2]
        \arrow[from=5-2, to=5-3]
        \arrow[from=3-1, to=2-1]
        \arrow[from=5-3, to=3-3]
        \arrow[from=3-2, to=3-3]
    \end{tikzcd}\]
\end{definition}

\begin{example}
    Locales are -- perhaps -- the best example of quantales that 
    are commutative, idempotent, integral (and hence both 
    semicartesian and right-sided), divisible and strong (both 
    trivially). Among which, and of special significance to 
    Sheaf Theory, are the locales of open subsets of a topological 
    space \(X\), where the order relation is given by the
    inclusion, the supremum is the union, and the finitary infimum 
    is the intersection.
\end{example}

\begin{example}\label{example: unital quantales}
    We list below some examples of unital quantales that are not 
    locales:
    \begin{enumerate}
        \item
            The extended half-line \([0,\infty]\) with 
            order the inverse order -- \(\geq\) --, and the 
            usual sum of real numbers as the 
            multiplication. 

            Since the order relation is \(\geq\), the top 
            element is \(0\) and the bottom elements is 
            \(\infty\). We call this the Lawvere quantale
            due to its relation to Lawvere spaces (related to metric spaces);
        \item 
            The extended natural numbers 
            \(\mathbb{N}\cup \{\infty\}\), as a restriction
            of the the Lawvere quantale  (related to 
            distance on graphs);
        \item 
            The set \(\mathcal{I}(R)\) of ideals of a 
            commutative and unital ring \(R\) with order 
            \(\subseteq\), and the multiplication as the 
            multiplication of ideals. The supremum is the 
            sum of ideals, the top element is \(R\) and the 
            trivial ideal is the bottom;
        \item 
            The set \(\mathcal{RI}(R)\) of right (or left) 
            ideals of an unital ring \(R\) with the same 
            order and multiplication of the above example. 
            Then the supremum and the top and the bottom 
            elements are also the same of 
            \(\mathcal{I}(R)\);
        \item 
            The set of closed right (or left) ideals of a 
            unital \(C^*\)-algebra, the order is the 
            inclusion of closed right (or left) ideals, and 
            the multiplication is the topological closure 
            of the multiplication of the ideals.
    \end{enumerate}
    For more details and examples we recommend 
    \cite{rosenthal1990quantales}. 
\end{example}

\begin{remark}
    The first three examples we introduced in 
    \ref{example: unital quantales} are commutative and integral 
    quantales. The last two examples are neither commutative 
    nor semicartesian. The forth is not idempotent but the fifth 
    is, and both are right-sided (resp. left-sided) quantales 
    \cite[cf.]{rosenthal1990quantales}.

    The main examples of strong quantales are Heyting algebras, 
    \(([0,1],\leq,\cdot)\), strict linear quantales. 
    Some MV-Algebras, like the Chang's 
    \(([0,1],\wedge, \vee, \oplus, \tensor, 0, 1)\) \cite{chang}'s,
    are not strong.\\

\end{remark}

\begin{center}
    \textbf{
        From now on, we assume all quantales to be commutative.
    }
\end{center}

\begin{definition}
    Let \(\quantale Q\) be a a quantale, we define an alternative 
    partial order \(\preceq\) given by
    \[
        a\preceq b \iff a = a \tensor b
    \]
\end{definition}

\begin{remark}~\label{remark:semicartesian properties}
\begin{enumerate}
    \item 
        Let \(\quantale Q\) be a unital quantale, then it is 
        integral iff it  is semicartesian;
    \item 
        \(a\preceq b\leq c\) implies \(a\preceq c\);
    \item 
        If \(e \in \extent\quantale Q\), 
        \((e \preceq a) \iff (e \leq a)\). 
\end{enumerate}
\end{remark}

\begin{prop}\label{prop:idempotent + semicartesian = locale}
    If \(\quantale Q\) is semicartesian and idempotent, it is in 
    fact a complete distributive lattice and \(\tensor=\wedge\). 
    In other words, it is a locale.
\end{prop}
\begin{proof}
    Suppose \(\quantale Q\) is in fact idempotent, we have --
    because \(\tensor\) is increasing in both arguments -- that
    \[
        a\leq b \implies a\leq c \implies a\leq b\tensor c
    \]

    Hence, if \(a\) is less than both \(b\) and \(c\), then it 
    must be smaller than \(b\tensor c\); but 
    since \(\quantale Q\) is semicartesian, by remark 
    \ref{remark:semicartesian properties} above, 
    \(\tensor\leq\wedge\). This means that \(b\tensor c\) is 
    a lowerbound for \(b\) and \(c\), but what we had gotten from 
    idempotency means it's the greatest such upper bound.

    Thus the multiplication satisfies the universal property of 
    infima. The above is just a particular case of 
    \cite[Proposition 2.1]{nlab:quantale}.
\end{proof}
	\section{Preliminaries on \(\quantale Q\)-Sets} 
\begin{remark}
    Hereon we are working exclusively with commutative 
    semicartesian quantales, as opposed to ones without those
    properties.
\end{remark}

Given a quantale \(\quantale Q\), one may form -- roughly speaking -- a
\(\quantale Q\)-space, wherein distances between points are measured by 
elements of \(\quantale Q\) as opposed to -- say -- \([0,\infty)\) as 
we often do. This definition is made precise in the notion of a 
\(\quantale Q\)-set.
\begin{definition}
    A \(\quantale Q\)-set is a set endowed with a 
    \(\quantale Q\)-distance operation usually denoted by 
    \(\delta\). \lat{i.e} a \(\quantale Q\)-set is a set \(X\) with
    a map \(\delta:X^2\to X\) satisfying:
    \begin{enumerate}
        \item \(\delta(x,y) = \delta(y,x)\);
        \item \(\delta(x,y)\tensor\delta(y,z)\leq\delta(x,z)\);
        \item \(\delta(x,x)\tensor\delta(x,y) = \delta(x,y)\).
    \end{enumerate}
    and it is usual to denote \(\delta(x,x)\) by simply the 
    “extent of \(x\)” written as \(\extent x\).
\end{definition}

A couple of things might jump out to reader in the above definition. 
\begin{enumerate*}[label=(\roman*)]
    \item	\(\delta\) is symmetric, even though we have thrown out 
        all but the vaguest notions tying ourselves to metric 
        spaces;
    \item	Why is the triangle inequality upside down?
    \item	\(\extent x\tensor\delta(x,y) = \delta(x,y)\),
        why not just ask that \(\extent x = \top\)?
\end{enumerate*}

Those questions are all valid -- and answering the first and last ones 
differently has been done in the past and is the main difference 
between \(\quantale Q\)-sets and \(\quantale Q\)-enriched categories 
from a definitional perspective. The question of order being inverse is 
more one of sanity: since we treat a \(\quantale Q\)-set as a set with
\(\quantale Q\)-valued equality, it makes sense to think that 
\(\extent x\) is the maximally valid equality to \(x\) and hence the
triangular inequality needs to be turned upsidedown -- and turned into
the transitivity of equality.

\begin{remark}
    When we speak of properties of the type 
    \(P(\vec x)\leq Q(\vec x)\) in \(\quantale Q\)-sets, it is 
    often more insightful to think that the logically equivalent 
    (but notationally less helpful) statement 
    \[
        P(\vec x) \residue Q(\vec x) \ (=\top)
    \]
\end{remark}

There are two main category structures that one can canonically endow
the collection of all \(\quantale Q\)-sets with. One is taking maps to
be co-contractions (\lat{i.e.} they make \(\delta\) bigger) -- the 
other is to consider well behaved \(\quantale Q\)-valued relations 
between the underlying sets.

\begin{definition}
    A functional morphism \(f:X\to Y\) is a function \(f\) between
    the underlying sets of \(X\) and \(Y\) such that \(f\) 
    increases \(\delta\) and preserves \(\extent\); that is to say
    \[
        \delta_X \leq \delta_Y\circ (f\times f) 
    \]
    \[
        \extent_X = \extent_Y \circ f
    \]

    A relational morphism \(\phi:X\to Y\) is a function 
    \(\phi:|X|\times|Y|\to\quantale Q\) satisfying 

    \begin{align*}
        \delta(x,x')\tensor\phi(x,y)
    &\leq	\phi(x',y)
    \\
        \phi(x,y) \tensor \delta(y,y') 
    &\leq	\phi(x,y')
    \\
        \phi(x,y) \tensor \phi(x,y') 
    &\leq	\delta(y,y')
    \\
        \phi(x,y) \tensor \extent y 
    &=	\phi(x,y)
    \\
        \extent x\tensor \phi(x,y)
    &=	\phi(x,y) 
    \\
        \bigvee_{y \in Y} \phi(x,y) 
    &=	Ex
    \end{align*}

    The reader should beware we don't often distinguish between 
    \(\delta_X\) and \(\delta_Y\) and instead rely on suggestively 
    named variables so as to indicate their type and hence the 
    \(\delta\) they refer to. In other words, the reader is 
    expected to be familiar with Koenig lookup%
    \footnote{
        ~Which, to quote a great website -- 
        \url{cppreference.com} -- 
        \begin{quote}
            Argument-dependent lookup, also known as ADL, 
            or Koenig lookup, is the set of rules for 
            looking up the unqualified function names in 
            function-call expressions, including implicit 
            function calls to overloaded operators.
        \end{quote}
    }.

    We denote by \(\QSet_r\) the category \(\quantale Q\)-sets 
    and relational morphisms between them and by \(\QSet_f\) the 
    category with the same objects but functional morphisms between
    them instead.
\end{definition}
\begin{prop}
    It should be observed that both notions of morphism actually do
    form a category. For functional morphisms, we take composition 
    to be the usual function composition -- since functional 
    morphisms are obviously closed under it -- and the identity 
    function as the idendity. For relational morphisms, the 
    identity on \(X\) becomes \(\delta_X\) and composition is the 
    (perhaps obvious) relational composition:
    \[
        [\psi\circ\phi](x,z) = \bigvee_{y\in Y}
            \phi(x,y)\tensor\psi(y,z)
    \]

    Since proving that functional morphisms do indeed form a 
    category would be trivial, we shall instead prove a stronger 
    result afterwards.
\end{prop}
\begin{proof}
    Firstly we should observe that \(\delta\) is indeed a 
    relational morphism. This is easy enough, as the first three 
    axioms are direct applications of the triangular inequality, 
    the fourth and fifth are direct applications of the 
    extension axiom of \(\quantale Q\)-sets, and the last one is 
    trivially true once we realize that
    \[
        \delta(x,y)\leq\extent x = \delta(x,x)
    \]

    Once we have that \(\delta\) is indeed a morphism, we can 
    wonder about the action of (pre-)composing with it. In which 
    case we obtain
    \begin{align*}
        [\delta\circ\phi](x,y) 
        &=	\bigvee_{y'\in Y}
            \phi(x,y')\tensor\delta(y',y)
        \\&\geq	\phi(x,y)\tensor\delta(y,y) 
        \\&=	\phi(x,y) 
        \intertext{
            on the other hand, we could have applied the 
            appropriate relational morphism axiom inside 
            the \(\bigvee\) thus getting that the composite 
            was smaller than \(\phi(x,y)\) instead. This 
            proves they were in fact equal all along. The 
            same goes for the other composite:
        }
        [\phi\circ\delta](x,y)
        &=	\bigvee_{x'\in X}
            \delta(x,x')\tensor\phi(x',y)
        \\&\geq	\delta(x,x)\tensor\phi(x,y)
        \\&=	\phi(x,y)
    \end{align*}

    Hence it only remains to see that \(\circ\) as defined for 
    relational morphisms is indeed a composition, in that it is 
    associative and that the composite of two morphisms is a 
    morphism. Firstly suppose that \(\phi\) and \(\psi\) are 
    composable morphisms -- we ought to show that their composite 
    is also a morphism:
    \begin{align*}
        [\psi\circ\phi](x,z)\tensor\delta(z,z')
        &=	\bigvee_{y}\phi(x,y)\tensor\psi(y,z)
                \tensor\delta(z,z')
        \\&\leq	\bigvee_{y}\phi(x,y)\tensor\psi(y,z')
        \\&=	[\psi\circ\phi](x,z')
        \intertext{
            And likewise, \lat{mutatis mutandis} one can 
            show that axioms 1-2 and 4-5 hold. Axiom 6 can
            be seen to hold easily as well:
        }
        \bigvee_z[\psi\circ\phi](x,z) 
        &=	\bigvee_y\phi(x,y)\tensor\bigvee_z\psi(y,z)
        \\&=	\bigvee_y\phi(x,y)\tensor\extent y
        \\&=	\bigvee_y\phi(x,y)
        \\&=	\extent x
        \intertext{
            Axiom 3 actually requires commutativity, which
            is unfortunately the first of many times it is 
            very necessary to make use of it;
        }
        [\psi\circ\phi](x, z)\tensor
        [\psi\circ\phi](x,z')
        &=	\bigvee_y\bigvee_{y'}
            \phi(x,y) \tensor\psi(y,z)\tensor
            \phi(x,y')\tensor\psi(y',z)
        \\&=	\bigvee_y\bigvee_{y'}
            \psi(y,z) \tensor
            \phi(x,y) \tensor
            \phi(x,y')\tensor
            \psi(y',z)
        \\&\leq	\bigvee_y\bigvee_{y'}
            \psi(y,z)   \tensor
            \delta(y,y')\tensor
            \psi(y',z)
        \\&\leq	\bigvee_y\bigvee_{y'}
            \psi(y,z)   \tensor
            \psi(y,z)
        \\&\leq	\bigvee_y\bigvee_{y'}
            \delta(z,z') = \delta(z,z')
    \end{align*}
    Associativity can obviously be seen to hold when we realize 
    that one can rearrange the terms of 
    \((\psi\circ\phi)\circ\chi\) and \(\psi\circ(\phi\circ\chi)\)
    to be in the form of 
    \[
        \bigvee_{x,y}
            \chi(w,x)\tensor\phi(x,y)\tensor\psi(y,z)
    \]
\end{proof}

\begin{definition}
    Instead of proving the category axioms for functional morphisms
    we promised to prove a stronger result -- which is incidently 
    useful for another paper of ours in the works -- which is to 
    prove that \(e\)-morphisms form a category (given a generic
    commutative unital quantale) and that functional morphisms form
    a wide subcategory of \(I\)-morphisms.

    So, let us define \(e\)-morphisms: given an idempotent element 
    \(e\) of \(\quantale Q\), a \(e\)-morphism is a 
    functional morphism “up to error \(e\)”:
    \[
        e\tensor\delta(x,x') \leq \delta(f(x),f(x'))
    \]
    \[
        \extent f(x) = e\tensor\extent f(x)
    \]
\end{definition}
\begin{prop}
    We claim that the collection of \(\langle e,\phi\rangle\) where
    \(\phi\) is an \(e\)-morphism constitutes a category under the 
    obvious composition laws. Furthermore, the identity function is 
    a \(I\)-morphism where \(I\) is the unit of the quantale, and 
    further still: \(I\)-morphisms are closed under composition and
    form a subcategory which is definitionally equal to 
    \(\QSet_f\).
\end{prop}
\begin{proof}
    Firstly, the obvious composition takes an \(e\)-morphism \(f\)
    and an \(e'\)-morphism \(g\) to a \((e\tensor e')\) morphism
    \(g\circ f\). Associativity is due to functional (in \(\Set\), 
    that is) \(\circ\) associativity and the fact that \(\tensor\)
    makes \(\quantale Q\) a semigroup. The fact that \(g\circ f\) 
    is a \((e\tensor e')\)-morphism is rather obvious and the proof
    is ommited.
    
    The identity is evidently a \(I\)-morphism -- and of course 
    that composing \(I\)-morphisms gives a 
    \(I\tensor I = I\)-morphism.
\end{proof}

\subsection{Some Examples}

\begin{example}[Initial Object]
    The empty set is -- vacuously -- a \(\quantale Q\)-set, and since
    morphisms are functions, it also happens to be the initial object.
\end{example}

\begin{example}[Terminal Object]
    The set of idempotent elements of \(\quantale Q\), denoted 
    \(\extent\quantale Q\) naturally has a structure of a
    \(\quantale Q\)-set -- given by 
    \(\delta(a,b)=a\tensor b=a\wedge b\)\footnote{This is the infimum in the subposet \(\extent\quantale Q\)m that turns out to be  a locale.}. It is trivial to see 
    that \(\tensor\) satisfies all \(\quantale Q\)-set laws. More 
    interestingly, however, \(\extent\quantale Q\) must be the 
    terminal object because \(\extent e = e\): since functional 
    morphisms preserve extents, one has that 
    \(\extent f(x) = \extent x\) -- however, \(\extent f(x) = f(x)\) 
    and thus \(f(x) = \extent x\). This proves that there is at most
    one morphism \(X\to\extent\quantale Q\).

    On the other hand, \(\delta(x,y)\leq\extent x\wedge\extent y\) 
    which just happens to make \(x\mapsto\extent x\) a functional 
    morphism. And thus, \(\extent\quantale Q\) is a \(\quantale Q\)%
    -set.
\end{example}

\begin{remark}
    One cannot use the whole of \(\quantale Q\) in 
    \(\extent\quantale Q\)'s stead, as 
    \(\delta(x,y)\tensor\extent x\) would not hold. The only reason 
    it holds in the above example is because for idempotents 
    \(\tensor=\wedge\). However, one can obtain a \(\quantale Q\)-set
    where the underlying set is \(\quantale Q\) itself
\end{remark}

\begin{example}
    Much akin to how Lawvere's quantale is a Lawvere space, a 
    (integral and commutative) quantale \(\quantale Q\) is a 
    \(\quantale Q\)-set. This is achieved with the following:
    \[
        \delta(x,y) = (x\residue y)\wedge(y\residue x)
    \]
    which is roughly equivalent to \(|x - y|\) for real numbers.
    This isn't necessarily the best \(\quantale Q\)-set structure
    we can give them, as \(\extent x = \top\) for any \(x\). 

    Ways to mitigate this phenomenon, which is specially strange 
    for \(\bot\), involve taking into account idempotents above 
    \(x\). An important quantallic property is the existence of 
    an operation \((\blank)^-\) taking an element \(x\) to the value 
    \(\sup\set{e\in\extent\quantale Q}{e\preceq x}\). Multiplying 
    \(\delta(x,y)\) by \(x^-\tensor y^-\) guarantees -- for instance
    -- that the above construction coincides with the terminal 
    object when \(\quantale Q\) is a locale.

    Another way to correct this, is to incorporate 
    \(\extent\quantale Q\) more directly, considering the space
    with underlying set \(\quantale Q \times\extent\quantale Q\)
    and \(\delta\) given by
    \[
        \delta((x,e),(y,a)) = 
        a\tensor e\tensor{[(x\residue y)\wedge(y\residue x)]}
    \]
    We write this \(\quantale Q\)-set as \({\quantale Q}_{\extent}\).
\end{example}

\begin{example}
    A construction that is explored in this work's sister-article 
    \cite{AMM2023constructions} but deserves to be mentioned here
    in passing \(X\chaosor X\), given by the underlying set 
    \(|X|\times|X|\) and with \(\delta\) given by the product of 
    the \(\delta\) of the coordinates. The reason this construction
    is relevant here is because 
    \(\delta: X\chaosor X\to{\quantale Q}_{\extent}\)
    in a natural way:
    \[
        (x,y)\mapsto (\delta(x,y), \extent x\tensor\extent y)
    \]
    And this happens to be a functional morphism.
\end{example}

\begin{remark}
    Monomorphisms are always injective functions, and epimorphisms 
    are always surjective. 
\end{remark}

\begin{example}[Regular Subobjects]
    A monomorphism is regular when it is an equalizer of a pair of 
    parallel arrows. Suppose \(f,g:A\to B\), one way to conceive of 
    an equalizer is as being the maximal subobject of \(A\) making 
    the diagram commute. It is quite trivial to see that subobjects
    in general \emph{must} be subsets with a point-wise smaller 
    \(\delta\). Hence, the largest subobject of \(A\) that 
    equalizes the pair is simply the subset of \(A\) where the 
    functions agree on, with the largest possible \(\delta\): that 
    being \(A\)'s. 

    Hence, we see a pattern where equalizers -- in general -- are 
    simply subsets with \(\delta\) coming from a restriction of 
    the original object. Importantly, though, we refer to 
    “regular subobjects” as monomorphisms that preserve \(\delta\)
    -- as they have been equivalently characterized.

    The skeptical reader might question that we have merely shown
    that regular monos preserve \(\delta\), as opposed to showing 
    such to be a sufficient condition. 
    
    In which case, given the fact that monos are injective 
    functions, \(\delta\)-preserving monos are simply subsets with 
    its superset's \(\delta\); consider one such mono: 
    \(m:A\hookrightarrow X\), with \(A\subseteq X\), one takes 
    \((X\amalg X)/\sim\) with \(\sim\) defined so as to identify 
    both copies of \(A\). 
    \[
        \delta(\eqclass{(x,i)},\eqclass{(y,j)}) = 
        \begin{cases}
            \delta(x,y), & i = j\\
            \bigvee_{a\in A}\delta(x,a)\tensor\delta(a,y), 
                & i\not = j
        \end{cases}
    \]
    There are two obvious inclusions of \(X\) into this set, 
    namely the upper branch and the lower branch. And they 
    coincide exactly on the section corresponding to \(A\), 
    so the equalizer of those two arrows must be the subset 
    inclusion of \(A\) into \(X\), the \(\delta\) must be 
    the biggest possible, which is the restriction of \(X\)'s
    so as to remain a morphism.
\end{example}




\begin{example}
     Suppose \((X,d)\) is a pseudo-metric space, 
     then \((X,d)\) is a \([0,\infty]\)-set where 
     \(\delta(x,y) \neq \bot = \infty, \forall x,y \in X\).
\end{example}


\begin{example}
    Given a nonempty index set \(I\), we have a 
    \(\quantale Q\)-set \(\coprod_{i\in I}\top\), given by
    \[
        \delta((e, i), (e', i')) = \begin{cases}
            e\wedge e', &\text{if} ~ i = i';\\
            \bot,       &\text{otherwise}.
        \end{cases}
    \]    
\end{example}

\begin{example}
    Given a commutative ring \(A\), let the set of its (left) ideals 
    be denoted \newcommand\ideals{\mathscr I} \(\ideals_A\). 
    \(\ideals_A\) is a quantale. Given a left \(A\)-module \(M\), we 
    can endow it with the structure of a \(\ideals_A\)-set:
    \[
        \delta(x,y) = \bigvee\set{
            I \in\ideals_{A}
        }{
            I\cdot x = I\cdot y
        }
    \]
    In fact, that supremum is attained at with a particular ideal. 
    Moreover, \(Ex = A = \max\ideals_A\). 
\end{example}
\begin{example}
    Suppose that \(\quantale Q\) is a quantale with 
    “idempotent upper  approximations”\footnote{%
        In \cite{tenorio2022sheaves} are described sufficient conditions
        for \(\quantale Q\) to have such a property.
    }: 
    \[
        \forall q\in\quantale Q:\exists q^+\in\extent\quantale Q:
            q^+ = \min\set{e\in\extent\quantale Q}{q\preceq e}
    \]
    Then
    \[
        \delta(x, y) = \begin{cases}
            x\tensor y, & x\not=y;\\
            x^+,        & x  =  y.
        \end{cases}
    \]
    defines a \(\quantale Q\)-set structure on \(\quantale Q\) itself.
\end{example}


\subsection{The Underlying Graph Functor}

\begin{definition}
    There is a functor \(\Graph:\QSet_f\to\QSet_r\) which is the 
    underlying \emph{graph} functor, because it takes functional 
    morphisms to their graph relation. More precisely, 
    \(\Graph(X) = X\) and given \(f:X\to Y\)
    \[
        (\Graph f)(x,y) = \delta(f(x), y) 
    \]
\end{definition}
\begin{prop}
    As defined above, \(\Graph f\) is indeed a functional morphism
    and \(\Graph\) is indeed a functor.
\end{prop}
\begin{proof}
    It is clear, from \(f\) being a functional morphism, that 
    \(\Graph f\) should at least satisfy the \(\delta\) and 
    \(\extent\) axioms for relational morphism axioms; the 
    (\(\Sigma\)) axiom holds because of the triangular inequality; 
    the strictness axiom holds as taking \(b = f(a)\) gives 
    \(\extent f(a)\) which is \(\extent a\).

    Regarding functoriality, \(\Graph\id = \delta\) and hence the
    identity in the relational category; moreover, 
    \begin{align*}
        (\Graph g)\circ(\Graph f)(x,z)
        &=	\bigvee_y
                (\Graph f)(x,y)
                \tensor
                (\Graph g)(y,z)
        \\&=	\bigvee_y
                \delta(f(x), y)
                \tensor
                \delta(g(y), z)
        \\&=	\delta(g\circ f(x), z)
        \\&=	\Graph(g\circ f)
    \end{align*}
\end{proof}

    \section*{Interlude}
	In the general course of studying \(\QSet[\Omega]\) for a frame/locale
	\(\Omega\) one is bound to eventually come across a certain notion 
	called “completeness”. Completeness is, largely speaking, a property 
	necessary to transition between \(\QSet[\Omega]\)-world to 
	\(\Sh\Omega\)-world 
	\cite[prop.~2.7.10,~cf. def. 2.9.1]{borceux1994handbook3}. This notion
	of completeness we call Scott-completeness and refers to a thing called
	singletons.

	Sheaves, almost by definition, are about gluing partial data -- in that
	morally speaking, compatible elements over a covering admit exactly one
	gluing which lies over the covered open. This gives a second natural 
	notion of completeness -- which we have come to call 
	Gluing-completeness. 

	In these next sections, therefore, we shall give the precise definitions 
	pertaining to both notions, make explicit some of the relations between
	them.

	\section{Gluing Completeness}
\begin{definition}[Compatible Family]
    Given a \(\quantale Q\)-set \(X\), a family \(A\) of its 
    elements is said to be \emph{compatible} when 
    \(\forall a, b\in A\),
    \[
        \delta(a, b) = \extent a \tensor \extent b
    \]
    naturally, since extents are idempontent, and for idempontents 
    \(\tensor=\wedge\), this is the same as saying 
    \(\delta(a,b) = \extent a\wedge\extent b\). 
\end{definition}

Note that the empty family of a \(\quantale Q\)-set is vacuously compatible.

In light of our usual interpretation for \(\delta\) being ‘the place 
where two fragments of data agree’, this says that \(a\) and \(b\) 
agree in the interesction of their extents -- their domains of 
definition and degrees of certainty, so to speak. Hence, they are 
compatible.
\begin{definition}[Gluing Element]
    Given a \(\quantale Q\)-set \(X\) and a family of its elements
    \(A\), we say that \(x_A\in X\) is a \emph{gluing of \(A\)} or 
    that \emph{it glues \(A\)} when \(\forall a\in A\),
    \[
        \delta(a, x_A) = \extent a
    \]
    \[
        \extent x_A = \bigvee_{a\in A}\extent a
    \]
\end{definition}
This definition then captures the idea that an amalgamation of local 
data agrees integrally with the data it agregates, and -- importantly 
-- does not make choices or provide information not present in the 
family \(A\) it glues.

\begin{prop}
    If a family \(A\) admits a gluing element, it is compatible.
    Moreover, any gluing elements are “morally unique” -- that is
    to say, they are \(\delta\)-equivalent.
\end{prop}
\begin{proof}
    Suppose \(x_A\) glues \(A\); then let \(a\) and \(b\) be 
    members of the family \(A\); observe that
    \begin{align*}
        \delta(x_A, a) &= \extent a\\
        \delta(x_A, b) &= \extent b
        \intertext{hence}
        \extent a \tensor\extent b &= 
        \delta(x_A, a)\tensor\delta(x_A, b) 
        \intertext{but then,}
        \extent a\tensor\extent b &\leq \delta(a,b)
        \intertext{but we know the following always holds}
        \delta(a,b)&\leq\extent a\tensor\extent b
        \intertext{therefore they are the same.}
    \end{align*}
    This suffices for compatibility, now for the second claim:
    consider \(y\in X\) and let us compare two possible gluings
    -- say \(x\) and \(x'\).
    \begin{align*}
        \delta(x, y) 
          &=	\delta(x,y)\tensor\extent x
        \\&=	\bigvee_{a\in A}\delta(x,y)\tensor\extent a
        \\&=	\bigvee_{a\in A}\delta(x,y)\tensor\extent a
                                       \tensor\extent a
        \\&=	\bigvee_{a\in A}\delta(x,y)\tensor\delta(x',a)
                                       \tensor\delta(x, a)
        \\&\leq	\bigvee_{a\in A}\delta(x,y)\tensor\delta(x',x)
        \\&\leq	\delta(x',y)
    \end{align*}
    since the same holds swapping \(x\) for \(x'\) and vice versa,
    it follows that \(\delta(x,y) = \delta(x',y)\) for \emph{any}
    \(y\); hence, they are essentially the same.
\end{proof}

\begin{definition}[\(\delta\)-equivalence]
    The above property is somewhat important, as -- as often is the
    case -- equality is \emph{evil} and equivalence is \emph{good}.
    Two elements \(x, y\) of a \(\quantale Q\)-set \(X\) are said 
    to be \emph{\(\delta\)-equivalent} whenever either of the 
    following equivalent conditions hold:
    \[
        \forall z\in X: \delta(x,z) = \delta(y,z)
    \]
    \[
        \delta(x,y) = \extent x = \extent y
    \]

    That the former implies the latter is rather obvious, simply 
    substitute \(z\) for \(y\) and for \(x\). That the latter 
    yields the former is also simple to show:
    \begin{align*}
        \delta(x,z) 
          &=	\delta(x,z)\tensor\extent x 
        \\&=	\delta(x,z)\tensor\delta(x,y)
        \\&\leq	\delta(y,z)
    \end{align*}
    and \lat{mutatis mutandis} one obtains the other required 
    inequality.
\end{definition}
\begin{remark}
    It is evident that this is an equivalence relation.
\end{remark}
Two \(\delta\)-equivalent points cannot -- for instance -- be taken
to non-equivalent points via a functional morphism; and relational 
morphisms can't distinguish them. This is codified in the following
lemma
\begin{lemma}
    The relation of \(\delta\)-equivalence is ‘congruential’ 
    for either type of morphism. In more precise terms: provided 
    with equivalent elements \(x\) and \(x'\); 
    \begin{enumerate}
        \item For functional morphisms, congruential means 
            \(f(x)\sim f(x')\); 
        \item For relational morphisms, it means that both
            \(\phi(x, y) = \phi(x',y)\) and 
            \(\psi(w,x) = \psi(w,x')\).
    \end{enumerate}
\end{lemma}
\begin{proof}
    \begin{align*}
        \extent x = \extent x' = \delta(x,x')
          &\leq	\delta(f(x),f(x'))
        \\&\leq	\extent f(x) = \extent x
        \\&=	\extent x' = \extent f(x')
        \intertext{As for relational morphisms,}
        \phi(x,y)
          &=	\phi(x,y)\tensor\extent x 
        \\&=	\phi(x,y)\tensor\delta(x,x')
        \\&\leq	\phi(x'y)
    \end{align*}
    The dual inequality holds for the same reasons; and similarly, 
    the same equality can be obtained for \(\psi\).
\end{proof}

We can state the gluing completeness as follows
\begin{definition}[Gluing Completeness]
    A \(\quantale Q\)-set \(X\) is said to be 
    \emph{gluing complete} when every compatible family has exactly
    one gluing. This is equivalent to asking that all compatible 
    families admit \emph{a} gluing and for 
    \(\delta\)-equivalence to be “extensionally equal to equality” 
    -- which sounds very pretentious.

    The equivalence can be seen to be true by realizing that 
    \(\{x,x'\}\) is glued by both \(x\) and \(x'\) if and only if 
    they are \(\delta\)-equivalent.
\end{definition}

\begin{definition}[Extensionality]
    Extensionality is the name of the property of 
    \(\delta\)-equivalence being extensionally equal to equality.
    That is, to say:
    \[
        \delta(x,y) = \extent x = \extent y \implies
        x = y
    \]
    or equivalently 
    \[
        \delta(\blank, x) = \delta(\blank, y)\implies
        x = y
    \]
\end{definition}

\begin{example}
    It is quite easy to see that the terminal \(\top\) 
    \(\quantale Q\)-set is gluing complete. Given that 
    \(\extent a \tensor\extent b = \delta(a,b)\) always 
    holds, any set of idempotent elements is compatible.
    Moreover, the gluing of such a set is its supremum, 
    as can easily be seen: 
    \[
        x_A = \extent x_A = 
        \bigvee_{a\in A}\extent a = 
        \bigvee_{a\in A} a
    \]
\end{example}

\begin{example}
    The initial object, the empty set, is trivially extensional
    but it fails to satisfy the gluing condition as the empty 
    family is also trivially \emph{compatible} and there is no 
    element in the empty set that actually glues it. Hence, 
    completeness doesn't hold, whereas extensionaly does.
\end{example}

\begin{example}
    In an extended pseudometric space \((X,d)\) seen as a 
    \([0,\infty]\)-set, points at infinity (with minimal extent)
    are always compatible with every other point, and normal 
    points are compatible if and only if they are morally the 
    same (ie. their distance is \(0\)). 
    
    Hence, one such extended pseudometric space satisfies the
    gluing condition if and only if it has a point at infinity,
    and is extensional if and only if it is morally a metric
    space (ie. it is a metric space if considered sans points
    at infinity) and it has exactly one point at infinity.
\end{example}

\begin{example}
    For nontrivial index sets \(I\) (namely: with more than one 
    element) and \(\quantale Q\)-sets \(X_i\) with elements
    \(\bot_i\) at infinity (null-extensional) 
    \(\coprod_{i\in I} X\) is never extensional, as 
    \(\delta((\bot_i, i),(\bot_j,j)) = \bot\) for multiple 
    different \(\bot_i\). In general, any \(\quantale Q\)-set
    with more than exactly one element with null extent cannot 
    be extensional. 
\end{example}


\begin{example} \label{deltaequiv-ex}
    \begin{enumerate}
    
    \item If \(A\) is a commutative ring and \(M\) is a left \(A\)-module, then the  set over \(\quantale Q = Ideals(A)\), \(X_M = (M, \delta)\), is such that for each \(x, y \in M\), \(x\sim_\delta y \)  iff \(A.x = A.y\), thus  it is not an extensional \(\quantale Q\)-set, in general: any member of a non-empty compatible family in \((M, \delta)\) is a gluing for that family. Moreover, \(M\) is a divisible \(A\)-module iff   \(\sim_\delta = \{(0,0)\} \cup (M\setminus\{0\})\times (M\setminus\{0\})\).

       \item If \(card(I) > 1\), then the \(\quantale Q\)-set \(I.1_{\quantale Q}-set\) is not extensional, since \(\delta((\bot,i),(\bot,j)) = \bot\) and there are \(i,j \in I\) such that \(i \neq j\). If \( I' \subseteq I\), then the family \(\{(\bot,i) : i \in I'\}\) is compatible and, for any \(j \in I\), \((\bot,j)\) is a gluing for that family.

    \item If \(\quantale Q\) is a integral quantale with  upper idempotent approximation, then the associated  \(\quantale Q\)-set on \(\quantale Q\) is extensional. A family \(S\) on that \(\quantale Q\)-set is compatible iff for each \(x,y \in S\), \(x \tensor y = x^+ \tensor y^+\). 
    
    
    \end{enumerate}
\end{example}

\begin{prop}
    The functor \(\Graph\) is faithful exactly on those pairs of 
    objects where the codomain is extensional. This is to say: if 
    \(Y\) is extensional, then the following map is injective for
    all \(X\)
    \[
        \QSet_f(X,Y) \xrightarrow \Graph \QSet_r(X,Y)
    \]

    Moreover, extensional \(\quantale Q\)-sets are the largest 
    class of objects such that you can quantify universally over
    \(X\).
\end{prop}
\begin{proof}
    First let's prove the latter claim: let \(Y\) be any 
    \(\quantale Q\)-set; and suppose \(y\) and \(y'\) are 
    \(\delta\)-equivalent. Now let \(X = \{*\}\), setting 
    \(\extent*=\extent y\).

    There are two obvious maps \(X\to Y\), namely \(*\mapsto y\) 
    and \(*\mapsto y'\). Call them \(f\) and \(g\) respectively.
    By the equivalent definition of \(\delta\)-equivalence, \(y\)
    and \(y'\) have the same \(\delta\)s, so
    \[
        \Graph(*\mapsto y) = \delta(\blank, y) = 
        \delta(\blank, y') = \Graph(*\mapsto y')
    \]

    And hence, if \(y \not= y'\), then \(\Graph\) must send 
    different functional morphisms to the same relational one
    and therefore fail injectivity.

    Now, that being said we still must prove that the functor is
    injective on \(\hom\)s when the codomain is extensional. But 
    this is simple. Suppose \(\Graph(f)=\Graph(g)\) -- and thus
    that for all \(x\) and \(y\), 
    \[
        \delta(f(x),y) = \delta(g(x), y)
    \]
    by taking \(y\) to be \(f(x)\) we get that 
    \(\delta(g(x),f(x))=\extent f(x)\) and \lat{mutatis mutandis}
    the same can be done to show that \(f(x)\) is always 
    \(\delta\)-equivalent to \(g(x)\). Because \(Y\) is extensional
    we have that \(g(x) = f(x)\) for all \(x\) and therefore they 
    are they one and the same functional morphism.
\end{proof}

\begin{remark}
    Completeness, be it gluing-wise or singleton-wise, is a rather 
    pointless notion in the category of relational morphisms; this
    is because every \(\QSet_r\) is equivalent to its full 
    subcategories whose objects are complete in either notion.

    Hence, any categorical property we \emph{could ever prove} 
    about those subcategories would also be true of the larger 
    relational category.
    
    In fact, completeness can -- in a way -- be understood as the 
    property one can impose on objects so as to make their 
    functional morphisms representable functor “the same” as their 
    relational morphisms representable functor. This will be 
    elaborated in a subsequent section.
\end{remark}

In the light of the remark above, we shall focus ourselves on the 
functional morphisms category, since talking about complete 
\(\quantale Q\)-sets under relational morphisms is just a roundabout 
way of talking about general \(\quantale Q\)-sets under those same 
morphisms.

\subsection{Some Results about Gluings and Compatible Families}

\begin{lemma}
    Family compatibility is preserved by direct images of 
    functional morphisms.
\end{lemma}
\begin{proof}
    Given \(A\) compatible for \(X\) and \(f:X\to Y\); take 
    any \(a,b\in A\); it suffices to show that 
    \(\delta(f(a),f(b)) = \extent f(a) \tensor\extent f(b)\).
    \begin{align*}
        \delta(a,b) 
          &\leq	\delta(f(a),f(b))
        \\&=	\delta(f(a),f(b))\tensor\extent f(a)
        \\&=	\delta(f(a),f(b))\tensor\extent f(a)
                             \tensor\extent f(b)
        \\&=	\delta(f(a),f(b))\tensor\delta(a,b)
        \\&\leq	\delta(a,b)
        \tag{semi-cartesian}
    \end{align*}
\end{proof}
\begin{lemma}
    Gluing elements are preserved by direct images of 
    functional morphisms
\end{lemma}
\begin{proof}
    Should \(x\) glue \(A\), we already know that \(f[A]\) is 
    compatible -- in light of the previous lemma and the fact
    that a family that admits a gluing must be compatible -- so 
    it merely suffices to show that \(f(x)\) glues \(f[A]\). 
    Take \(a\in A\),
    \begin{align*}
        \extent a
          &=	\delta(x,a)
        \\&\leq	\delta(f(x),f(a))
        \\&\leq \extent f(a) = \extent a
        \intertext{Thus establishing the first condition.}
        \extent f(x) 
          &= \extent x 
        \\&= \bigvee_{a\in A}\extent a
        \\&= \bigvee_{a\in A}\extent f(a)
        \\&= \bigvee_{\alpha\in f[a]}\extent \alpha
    \end{align*}
    Thus completing the proof.
\end{proof}

We shall see that there is a natural way of ‘making’ a 
\(\quantale Q\)-set gluing complete. This what we shall focus ourselves
on in the next subsection.

\subsection{Gluing-Completion}
Completion problems often come in the shape of 
“all \(X\) induce a \(Y\), but not all \(Y\) come from an \(X\) -- can
we make it so \(X\) and \(Y\) are in correspondence?” and the solution 
to those kind of problems can often be found in looking at the 
collection of all \(Y\)s or quotients of it. 

In that spirit, we shall give the collection of all compatible families
the structure of a \(\quantale Q\)-set; we shall embed the original 
\(\quantale Q\)-set in it; and then we “throw away” redundant points in
an appropriate quotient.

\begin{definition}
    Given a \(\quantale Q\)-set \(X\), we say \(\Gluings(X)\) is 
    the \(\quantale Q\)-set given by
    \begin{align*}
        |\Gluings(X)| &= \set[:]{A\subseteq|X|}{
            A~\text{is compatible}
        }\\
        \delta(A, B) &= 
            \bigvee_{\substack{a\in A\\b\in B}}
            \delta(a,b)
    \end{align*}
    Obviously, we claim that the above does indeed define a 
    \(\quantale Q\)-set.
\end{definition}
\begin{proof}
    \begin{align*}
        \delta(A,B)\tensor\delta(B,C) 
          &=	\bigvee_{a\in A}
            \bigvee_{b\in B}
            \bigvee_{\beta\in B}
            \bigvee_{c\in C}
                \delta(a, b)\tensor\delta(\beta, c)
        \\&=	\bigvee_{a\in A}
            \bigvee_{b\in B}
            \bigvee_{\beta\in B}
            \bigvee_{c\in C}
                \delta(a, b)\tensor
                \extent b\tensor
                \extent \beta\tensor 
                \delta(\beta, c)
        \\&=	\bigvee_{a\in A}
            \bigvee_{b\in B}
            \bigvee_{\beta\in B}
            \bigvee_{c\in C}
                \delta(a, b)\tensor
                \delta(b,\beta)\tensor
                \delta(\beta, c)
        \\&\leq	\bigvee_{a\in A}
            \bigvee_{c\in C}
                \delta(a, c)
        \\&=	\delta(A,C)
        \intertext{Which suffices for transitivity/triangular 
        inequality. It is obvious that the definition is 
        symmetric, so it only remains to show that extents 
        behave as they should.}
        \delta(A, A)\tensor\delta(A,B) 
          &=	\bigvee_{a, a'\in A}
            \bigvee_{\alpha\in A, b\in B}
                \delta(a, a')
                \tensor
                \delta(\alpha, b)
        \\&\geq	\bigvee_{\alpha\in A, b\in B}
            \delta(\alpha,\alpha)
            \tensor
            \delta(\alpha, b)
        \\&=	\delta(A, B)
        \intertext{and thus}
        \tag{semi-cartesian}
        \extent A \tensor\delta(A,B)&=\delta(A,B)
    \end{align*}
    Thus proving it is actually a \(\quantale Q\)-set as desired.
\end{proof}

\begin{definition}[The \(\quantale Q\)-set of Compatible Families]
    We extend \(\Gluings\)\footnote{
        Here it seems necessary to justify the use of a fraktur
        typeface; Capital \textbf G is fine because we are 
        taking spaces with all possible \textbf{G}luings. But 
        since they aren't extensional, our functor is still a 
        bit ‘broken’. You can call it a mneumonic, but it's a 
        visual pun.
    } from an object-map to a dignified functor, by taking 
    functional morphisms to their direct image; this has the 
    surreptitious implicit premise that direct images are indeed
    functional morphisms between the compatible families 
    \(\quantale Q\)-sets.
\end{definition}
\begin{proof}
    \begin{align*}
        \delta(A, B) 
          &=	\bigvee_{a, b}\delta(a,b)
        \\&\leq	\bigvee_{a, b}\delta(f(a), f(b))
        \\&=	\delta(f[A], f[B])
    \end{align*}
    For extents it's even more immediate
    \[
        \extent A =\bigvee_{a\in A}\extent a = 
        \bigvee_{a\in A} \extent f(a)=\extent f[A]
    \]
\end{proof}

As mentioned in the previous footnote, \(\Gluings\) doesn't take 
objects to gluing-complete \(\quantale Q\)-sets; think of it in the 
same vein of object as Cauchy sequences sans the quotient. We first 
show that every compatible family over \(\Gluings(X)\) has a gluing.

\begin{lemma}
    If \(\mathcal A\) is a compatible family over \(\Gluings(X)\),
    then \(\bigcup\mathcal A\) is compatible over \(X\).
\end{lemma}
\begin{proof}
    Take \(A, B\in\mathcal A\) and \(a\in A\) and \(b\in B\); it 
    suffices to show that 
    \(\delta(a,b)\geq\extent a\tensor\extent b\) since we know the 
    dual inequality always holds.
    \begin{align*}
        \intertext{The following is always true:}
        \delta(a,b)
        &\geq	\bigvee_{a'\in A}
            \bigvee_{b'\in B}
                \delta(a, a')\tensor
                \delta(b, b')\tensor
                \delta(a', b')
        \intertext{forcing us to accept the remaining argument}
        &=	\bigvee_{a'\in A}
            \bigvee_{b'\in B}
                \extent a\tensor\extent a'\tensor
                \extent b\tensor\extent b'\tensor
                \delta(a', b')
        \\&=	\extent a\tensor\extent b\tensor
            \bigvee_{a',b'}\delta(a',b')
        \\&=	\extent a\tensor\extent b\tensor
            \delta(A,B)
        \\&=	\extent a\tensor\extent b\tensor
            \extent A\tensor\extent B
        \intertext{since for extents \(\tensor=\wedge\), and 
        since \(\extent a\leq\extent A\) and likewise for 
        \(b\) and \(B\),}
        &=	\extent a\tensor\extent b
    \end{align*}
\end{proof}

\begin{prop}[Compatible Families of \(\Gluings\) have Gluings]
    We must give at least \emph{one} gluing, but in general we 
    should have many such gluings. That's a defect we shall correct
    by taking an appropriate quotient later.
\end{prop}
\begin{proof}
    Given \(\mathcal A\), we claim that \(\bigcup\mathcal A\) is 
    a canonical candidate for its gluing. We shall now verify that
    it indeed is up for the task. First, let \(A\in\mathcal A\) and
    denote by \(X\) the family given by \(\bigcup\mathcal A\); let 
    us proceed.
    \begin{align*}
        \delta(X, A) 
        &=	\bigvee_{x\in X, a\in A}\delta(x,a)
        \\&\geq	\bigvee_{x\in A, a\in A}\delta(x,a)
        \\&=	\delta(A, A) = \extent A
        \\&\geq	\delta(X, A)
        \intertext{Remaining for us simply to show that the 
        extent matches what we wanted.}
        \extent X 
        &=	\bigvee_{x\in X}\extent x
        \\&=	\bigvee_{A\in\mathcal A}
            \bigvee_{x\in A}
                \extent x
        \\&=	\bigvee_{A\in\mathcal A}
                \extent A
    \end{align*}
\end{proof}

The last piece of the puzzle of gluing completion is ensuring 
extensionality; and this is done by another functor, the 
\(\delta\)-quotient functor.

\begin{definition}[\(\delta\)-quotient]
    We define the \(\faktor{-}{\delta}\) functor naturally by 
    quotienting together \(\delta\)-equivalent elements. It 
    does indeed take \(\quantale Q\)-sets to \(\quantale Q\)-sets
    if we take \(\delta(\eqclass{x},\eqclass{y}) = \delta(x,y)\).
    It does indeed take morphisms to morphisms, due to \(\delta\)
    equivalence being congruential.
\end{definition}

\begin{lemma}
    A compatible family of \(\faktor{X}{\delta}\) has a gluing if
    and only if (any) section of that family has a gluing in \(X\).
    That is to say, if \(A\) is our compatible family and 
    \(\bar A\) is such that \(\coproj[\bar A] = A\), then 
    \(A\) has a gluing if and only if \(\bar A\) does. Moreover, if
    any section has a gluing, then every section is also glued by 
    it.
\end{lemma}
\begin{proof}
    Suppose that \(A\) has a gluing, call it \(\eqclass x\); by 
    definition,
    \[
        \extent x = 
        \extent\eqclass x = 
        \bigvee_{\eqclass a\in A}\extent\eqclass a = 
        \bigvee_{\eqclass a\in A}\extent a
    \]
    So, if \(\bar A\) is a section of \(A\), then
    \[
        \extent x = \bigvee_{a\in\bar A}\extent a
    \]
    \[
        \delta(x, a) = \delta(\eqclass x,\eqclass a) = 
        \extent\eqclass a = \extent a
    \]
    Thus, \(x\) glues any section of \(A\). Now suppose instead 
    that \(x\) glues some section \(\bar A\) of \(A\). Of course 
    \(\eqclass x\) must now glue \(A\):
    \[
        \extent \eqclass x = \extent x = 
        \bigvee_{a\in\bar A}\extent a = 
        \bigvee_{a\in\bar A}\extent\eqclass a = 
        \bigvee_{\eqclass a\in A}\extent\eqclass a = \extent A
    \]
    \[
        \delta(\eqclass a, \eqclass x) = \delta(a, x) = 
        \extent a = 
        \extent\eqclass a
    \]
\end{proof}
\begin{remark}
    Evidently, every compatible family over \(\faktor X\delta\) has
    a section that is a compatible family over \(X\).
\end{remark}

\begin{lemma}
    The functor \(\faktor -\delta\) makes \(\quantale Q\)-sets 
    extensional (in fact it's “the best” extensional 
    approximation).
\end{lemma}
\begin{proof}
    Suppose 
    \[
        \delta(\eqclass x, \eqclass y) = 
        \extent\eqclass x = 
        \extent\eqclass y
    \]
    Then \(\delta(x, y) = \extent x = \extent y\) and immediately
    we have \(x\sim y\) and as such, \(\eqclass x = \eqclass y\).

    The claim that it is “the best” at doing this is harder to 
    state and prove, and it could/should be a theorem.
\end{proof}

\begin{theorem}[it became a theorem]
    The functor \(\QSet_{f\mathrm{Ext}}\hookrightarrow\QSet_f\) has
    a left adjoint and it is \(\faktor-\delta\).
\end{theorem}
\begin{proof}
    As will become a trend throughout this article, we go about 
    looking for units and counits and proving the triangle/zig-zag 
    identites for those. 

    The counit is obvious, as applying \(\faktor-\delta\) to an 
    already extensional object yields an isomorphic object whose 
    elements are singleton sets containing exactly the elements of
    the original \(\quantale Q\)-set. So we simply unwrap them and
    that's the counit.

    The unit is also simple, it's just the quotient map \(\coproj\)
    -- taking elements to their \(\delta\)-equivalence class.

    In this case, there is no real point to showing the zig-zag
    identities hold; but the reader may if they feel so inclined.
\end{proof}

\begin{definition}[Gluing Completion]
    The gluing completion functor \(\GluingCompletion\) is simply
    \((\faktor-\delta)\circ\Gluings\).
\end{definition}

\begin{prop}
    The Gluing-completion of a \(\quantale Q\)-set is indeed 
    gluing-complete.
\end{prop}
\begin{proof}
    We know that every compatible family over 
    \(\GluingCompletion X\) has at least one gluing, but we also 
    know they are extensional; and thus he have that they are 
    gluing complete.
\end{proof}

\begin{definition}[\(\GluingQSet\)]
    Here we give the name \(\GluingQSet\) to the full
    subcategory of \(\QSet_f\) whose objects are gluing complete.
    We drop the \(f\) subscript as this notion isn't particularly
    useful for relational morphisms anyway, so a it's pointless 
    datum. In a later subsection the \(f\) subscript will reappear
    as will a \(r\) subscript, so that we may prove that 
    \[
        \GluingQSet_r\equiv\QSet_r
    \]

    In the light of the previous proposition, we are entitled to 
    say that the signature of \(\GluingCompletion\) is actually
    \[
        \GluingCompletion:\QSet_f\to\GluingQSet
    \]
\end{definition}

\begin{lemma}
    Equivalent compatible families (in the sense of the canonical
    \(\delta\) we have defined for the set of such families) have 
    exactly the same gluings -- and they are equivalent to a 
    singleton set if and only if its element is one of their 
    gluings.
\end{lemma}
\begin{proof}
    It suffices that we show \(x\) glues \(A\) iff. \(\{x\}\sim A\)
    -- as the \(\sim\) being an equivalence takes care of the rest 
    of the claim for us.

    So, let's suppose that \(x\) does indeed glue \(A\); it follows
    that
    \[
        \delta(\{x\}, A) = 
        \bigvee_{a\in A}
            \delta(x,a) = 
        \bigvee_{a\in A}
            \extent a = \extent A = \extent \{x\}
    \]

    Now instead suppose that \(\{x\}\sim A\); In particular 
    \(\extent x = \extent \{x\}\) must be \(\extent A\) and thus
    \(\bigvee_{a\in A}\extent a\). So that takes care of the last 
    condition; the first condition for gluing still remains. But 
    not for long:
    \begin{align*}
        \extent a 
        &=	\extent a \tensor \extent A 
        \\&=	\extent a \tensor \delta(\{x\}, A)
        \\&=	\bigvee_{a'\in A}
            \extent a\tensor\delta(x, a') 
        \\&=	\bigvee_{a'\in A}
            \extent a\tensor\extent a'\tensor\delta(x, a')
        \\&=	\bigvee_{a'\in A}
            \delta(a, a')\tensor\delta(x, a')
        \\&\leq	\delta(x, a)
    \end{align*}
\end{proof}

\begin{theorem}[The Gluing Completion Adjunction]
    It's not particularly useful to have a functor that makes 
    \(\quantale Q\)-sets gluing-complete if that completion isn't
    actually “universal” in that it is the reflector of the 
    subcategory we are interested in.
\end{theorem}
\begin{proof}
    Again, we proceed by showing a unit-counit pair satisfying the
    zig-zag identities. For (my) sanity, let us name 
    \(L=\GluingCompletion\) and \(R\) the fully faithful inclusion
    \(\GluingQSet\hookrightarrow\QSet_f\).
    
    Let's find \(\eta:\id\to R\circ L\); what it must do is take an
    element \(x\) of \(X\) to something in the completion of \(X\);
    the completion of \(X\) is comprised of equivalence classes of 
    compatible families, so it is enough that we find a compatible 
    family to assign \(x\) to. \(\{x\}\) is the natural candidate 
    as
    \[
        \delta(x, y) = \delta(\{x\},\{y\}) = 
        \delta(\eqclass{\{x\}},\eqclass{\{y\}})
    \]
    Of course, we have only defined components, and don't know if 
    they form a natural transformation together; naturality holds
    quite trivially though:
    \[\begin{tikzcd}
        X &&& {\GluingCompletion(X)} \\
        & x & {\eqclass{\{x\}}} \\
        & f(x) & ? \\
        Y &&& {\GluingCompletion(Y)}
        \arrow["\eta", from=1-1, to=1-4]
        \arrow["\eta", from=4-1, to=4-4]
        \arrow["f"{description}, from=1-1, to=4-1]
        \arrow["{\GluingCompletion(f)}"{description}, 
            from=1-4, to=4-4]
        \arrow[maps to, from=2-2, to=3-2]
        \arrow[maps to, from=2-3, to=3-3]
        \arrow[maps to, from=3-2, to=3-3]
        \arrow[maps to, from=2-2, to=2-3]
    \end{tikzcd}\]
    We know that the action of \(\faktor-\delta\) on an arrow is 	
    \[
        f \xmapsto{\faktor-\delta}
        (\eqclass x \xmapsto{\eqclass{f}} \eqclass{f(x)})
    \]
    we know this from regular old set theory, it's just how 
    quotients work. The action of \(\Gluings\) on morphisms is 
    \[
        f \xmapsto\Gluings
        (A\xmapsto{f_*} f[A])
    \]
    that is, obviously, just taking functions to direct images.
    And thus, we combine the two to get
    \[
        \eqclass{f_*}(\eqclass{\{x\}}) = 
        \eqclass{f[\{x\}]} = 
        \eqclass{\{f(x)\}}
    \]
    making the above square commute.

    We should now provide a counit \(\epsilon:L\circ R\to\id\). 
    One such function should take things from the completion of an 
    already gluing-complete \(\quantale Q\)-set back to the 
    original \(\quantale Q\)-set. The things \(\epsilon\) operates 
    over are of the form \(\eqclass{A}\) where \(A\) is a 
    compatible family over -- say -- \(X\).

    The natural candidate is to take \(A\) to the unique element 
    \(x_A\) in \(X\) (which was gluing-complete to begin with) that
    glues it. Considering the previous lemma's, it is evident that 
    this definition does indeed yield a \emph{map}, as any 
    equivalent family will be glued by the same unique element.
    Moreover, since \(X\) was complete to begin with, 
    \(\{x_A\}\in\eqclass A\) and hence the map will obviously be a 
    morphism, as 
    \[
        \delta(\eqclass A,\eqclass B) = 
        \delta(\eqclass{x_A},\eqclass{x_B}) = 
        \delta(x_A, x_B)
    \]
    
    Remaining now only the need to show \(\epsilon\) is natural 
    transformation \(L\circ R\to\id\); 
    \[\begin{tikzcd}
        {\GluingCompletion(X)} &&& X \\
        & {\eqclass{A}} & {x_A} \\
        & {\eqclass{f[A]}} & {?} \\
        {\GluingCompletion(Y)} &&& Y
        \arrow["{\GluingCompletion(f)}"{description}, from=1-1,
            to=4-1]
        \arrow["\epsilon", from=4-1, to=4-4]
        \arrow["\epsilon", from=1-1, to=1-4]
        \arrow["f"{description}, from=1-4, to=4-4]
        \arrow[maps to, from=2-2, to=3-2]
        \arrow[maps to, from=3-2, to=3-3]
        \arrow[maps to, from=2-2, to=2-3]
        \arrow[maps to, from=2-3, to=3-3]
    \end{tikzcd}\]
    since gluing elements are preserved by functional morphisms, 
    the square does indeed commute; hence, the family of maps is 
    actually natural transformation and we can proceed in showing
    that the selected transformations satisfy the zig-zag 
    identities. Recapping, 
    \begin{align*}
        x&\xmapsto{\eta_X}\eqclass{\{x\}}\\
        \eqclass{A}&\xmapsto{\epsilon_{X'}}x_{A}
        \intertext{and we aim to show that}
        (\epsilon L)\circ(L\eta) &= \id\\
        (R\epsilon)\circ(\eta R) &= \id
    \end{align*}
    Given the below, it should then become obvious that we do 
    indeed have witnesses to the adjunction we have set out to 
    show to hold.
    \begin{align*}
        L\eta : L
        &\xrightarrow{\phantom{(L\eta)_X}}
        L\circ R\circ L\\
            (X\in\QSet)\mapsto( 
                \eqclass{A}
                &\xmapsto{(L\eta)_X}
                \eqclass{\{\eqclass A\}}
            )\\
        \\
        \epsilon L: L\circ R\circ L
        &\xrightarrow{\phantom{(\epsilon R)_X}}
        L\\
            (X\in\QSet)\mapsto(
                \eqclass{\{\eqclass A\}}
                &\xmapsto{(\epsilon L)_X}
                \eqclass A
            )\\
        \\
        R\epsilon: R\circ L\circ R
        &\xrightarrow{\phantom{R\epsilon)_K\ }}
        R\\
            (K\in\GluingQSet)\mapsto(
                \eqclass{A}
                &\xmapsto{(R\epsilon)_K}
                x_A
            )\\
        \\
        \eta R: R
        &\xrightarrow{\phantom{\eta R)_K\ }}
        R\circ L\circ R\\
            (K\in\GluingQSet)\mapsto(
                x
                &\xmapsto{(\eta R)_K}
                \eqclass{\{x\}}
            )
    \end{align*}
\end{proof}

\begin{remark}
    The adjunction counit is very easily seen to be a natural 
    \emph{isomorphism} (because \(R\) was fully faithful), however 
    the unit isn't: Consider the 
    \(\quantale Q\)-set given naturally by 
    \[
        X:=\{\bot\}\amalg\{\bot\}
    \]
    It obviously won't be extensional, because of this,
    \(\GluingCompletion(X)\iso\{\bot\}\) and there is no possible
    \emph{functional} isomorphism between those, because functional
    isomorphisms are -- in particular -- also bijections. So the 
    adjunction above is not an adjoint equivalence. For relational 
    morphisms there are those kinds of isomorphisms, and so the 
    adjunction unit will be a natural isomorphism too.
\end{remark}
	\section{Singleton/Scott-Completeness}
As previously mentioned, there is another notion completeness; it 
relates to a concept called singleton, which we will talk about now.

A singleton on a \(\quantale Q\)-set \(X\), morally speaking, is a 
\(\quantale Q\)-valued characteristic function. Although it isn't 
\emph{per se} a construction in the category \(\QSet\) -- it does lead 
to categorically relevant results.
\begin{definition}[Singleton]
    A \emph{singleton} on \(X\) is a map between the underlying set
    and \(\quantale Q\): \(\sigma:|X|\to\quantale Q\). This map is
    expected to satisfy the following axioms
    \begin{enumerate}
        \item \(\sigma(x)\tensor\extent x = \sigma(x)\);
        \item \(\sigma(x)\tensor\delta(x,y) \leq\sigma(y)\);
        \item \(\sigma(x)\tensor\sigma(y) \leq \delta(x,y)\);
        \item \(\sigma(x)\tensor\bigvee_y\sigma(y)=\sigma(x)\);
    \end{enumerate}
    the first two being called the “subset axioms” the third being 
    the “singleton condition” and the last axiom being called the
    “strictness condition”.
\end{definition}

Intuitively, singletons represent a membership relation to a set -- 
this is governed by the subset axioms; the singleton condition then 
states that the simultaneous membership of two different points implies 
their similitude; the last axiom, strictness, is more technical in 
nature. 

It has been known in the literature for quite a few years 
\cite[p. 30]{HOHLE199115} -- but its meaning is somewhat obscure, in 
that it is a necessary condition for a singleton to be representable 
(to be defined next), but it doesn't have a neat interpretation. The 
best we have arrived at is that strictness gives us that the supremum 
\(\sup_y\sigma(y)\) “behaves like an extent” in that it is idempotent.

\subsection{Basics about Singletons}
\begin{definition}[Representable Singleton]
    A singleton \(\sigma\) is said to be representable 
    (or represented by an element \(x_\sigma\)) when there is one
    such element \(x_\sigma\) such that 
    \[
        \sigma(y) = \delta(y,x_\sigma)
    \]
\end{definition}

\begin{prop}
    If \(A\) is a compatible family, then \(\delta(\{\blank\}, A)\)
    is a singleton. In particular, \(\delta(\blank, x)\) is always 
    a singleton (by taking \(A=\{x\}\)).
\end{prop}
\begin{proof}
    The subset axioms are immediately satisfied by the fact that 
    \(\Gluings(X)\) is always a \(\quantale Q\)-set. The Singleton
    condition is trivial too:
    \[
        \sigma_A(x) \tensor\sigma_A(y)= 
        \delta(\{x\}, A)\tensor\delta(\{y\}, A) \leq
        \delta(\{x\},\{y\}) = \delta(x,y)
    \]
    obviously, it is also strict considering that
    \[
        \bigvee_{y\in X}\sigma_A(y) \geq 
        \bigvee_{a\in A}\sigma_A(a) = 
        \bigvee_{a\in A}\extent a = 
        \extent A \geq
        \bigvee_{y\in X}\sigma_A(y)
    \]
\end{proof}
\begin{remark} 
    Quite obviously, if the singleton above is representable, 
    then the representing element glues \(A\). And if an element
    glues \(A\), then it is bound to represent the singleton.
\end{remark}

\begin{prop}
    \label{prop: singletons are continuous over compatible families}
    If \(A\) is a compatible family, and \(\sigma\) is some 
    singleton over \(A\)'s \(\quantale Q\)-set, we have that 
    if \(x_A\) glues \(A\), then 
    \[
        \sigma(x_A) = \bigvee_{a\in A}\sigma(a)
    \]
\end{prop}
\begin{proof}
    \begin{align*}
        \sigma(a) 
             =&\sigma(a)\tensor\extent a
        \\   =&\sigma(a)\tensor\delta(a, x_A)
        \\\leq&\sigma(x_A)
        \\
        \\
        \sigma(x_A) 
             =&\sigma(x_A)\tensor\extent x_A
        \\   =&\sigma(x_A)\tensor\bigvee_{a\in A}\extent a
        \\   =&\sigma(x_A)\tensor\bigvee_{a\in A}\delta(a, x_A)
        \\   =&\bigvee_{a\in A}\sigma(x_A)\tensor\delta(a, x_A)
        \\\leq&\bigvee_{a\in A}\sigma(a)
    \end{align*}
\end{proof}

\begin{definition}[Scott-Completeness]
    We say that \(X\) is Scott-complete, or singleton-complete, 
    when every one of its singletons have exactly one representing
    element.
\end{definition}

\begin{example}
    The terminal object is Scott-complete. An indirect way of 
    seeing this is to skip through the paper to where we prove 
    that the inclusion of Scott-complete \(\quantale Q\)-sets 
    creates limits. A more direct answer comes from the fact that
    since \(\top\in\extent\quantale Q\) is the gluing of the 
    compatible family \(\extent\quantale Q\), thanks to 
    prop. 
    \ref{prop: singletons are continuous over compatible families},
    \(\sigma(\top) = \bigvee_{e}\sigma(e)\). And hence:
    \[
        \sigma(e) = 
        \sigma(e)\tensor e \leq
        \sigma(\top)\tensor e = 
        \sigma(\top)\tensor\delta(e,\top) \leq 
        \sigma(e)
    \]
    And therefore, \(\sigma(\top)\) represents \(\sigma\).
    Uniqueness is quite obvious: if for some \(x\),
    \(\forall e:\sigma(e) = e\tensor x\), in particular
    \(\sigma(\top) = \top\tensor x = x\) by integrality.
\end{example}


\begin{example} \label{scottcompl-ex}
\begin{enumerate}
\item If \(A\) is a commutative ring and \(M\) is a left \(A\)-module, then the  set over \(\quantale Q = Ideals(A)\), \(X_M = (M, \delta)\), is such that for each \(x, y \in M\), \(x\sim_\delta y \)  iff \(A.x = A.y\), so in general, is not a extensional \(\quantale Q\)-set.
\end{enumerate}
\end{example}

\begin{remark}
    Scott-completeness is invariant under isomorphisms in 
    \(\QSet_f\) -- but not in \(\QSet_r\), this will be explained 
    in a later subsection. 

    The first claim is easy to show, as functional isomorphisms are 
    simply bijections preserving \(\delta\).
\end{remark}

\begin{theorem}
    Scott-Completness implies Gluing-Completeness
\end{theorem}
\begin{proof}
    Take a compatible family \(A\), consider the singleton defined
    in the previous proposition; suppose \(x_A\) represents 
    \(\sigma_A\); thus,
    \[
        \bigvee_{a\in A}\delta(x_A,a) =
        \sigma_A(x_A) = \delta(x_A, x_A) = \extent x_A
    \]
    \[
        \delta(a, x_A) = 
        \sigma_A(a) = 
        \bigvee_{a'\in A}\delta(a,a') = 
        \extent a 
    \]
    and hence \(x_A\) glues \(A\). Now suppose that \(x_A\) glues 
    \(A\), then 
    \begin{align*}
        \delta(y, x_A) 
        &\geq	\delta(x_A, a)\tensor
            \delta(a, y)
        \\&\geq	\extent a\tensor\delta(a,y) = \delta(a,y)\\
        \\
        \delta(y, x_A)
        &\geq	\bigvee_{a\in A}
                \delta(y,a)
        \\&\geq	\bigvee_{a\in A}
                \delta(a, x_A)\tensor
                \delta(y,x_A)
        \\&=	\delta(y,x_A)\tensor
            \bigvee_{a\in A}
                \delta(a, x_A)
        \\&=
            \delta(y,x_A)\tensor
            \bigvee_{a\in A}
                \extent a
        \\&=	\delta(y,x_A)\tensor\extent x_A
        \\&=	\delta(y,x_A)
    \end{align*}
    hence \(x_A\) represents \(\sigma_A\). Thus, if all singletons 
    have exactly one representing element, then all compatible 
    families have exactly one gluing.
\end{proof}

\begin{theorem}
    There are gluing complete \(\quantale Q\)-sets over 
    \(\quantale Q = \parts(2)\) which \emph{are not} Scott-complete
    -- so those two concepts \emph{are not} logically equivalent.
\end{theorem}
\begin{proof}
    Let \(\quantale Q\) be the following partial order
    
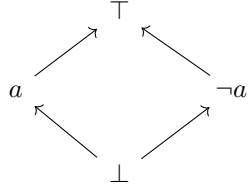
\[\begin{tikzcd}
        & \top \\
        a && {\neg a} \\
        & \bot
        \arrow[from=2-1, to=1-2]
        \arrow[from=2-3, to=1-2]
        \arrow[from=3-2, to=2-1]
        \arrow[from=3-2, to=2-3]
    \end{tikzcd}\]
    First, let \(S\subseteq\quantale Q\) and let \(\delta=\wedge\);
    \(S\) is gluing-complete if and only if \(S\) is complete as a
    lattice:

    If \(A\) is a subset of \(S\), it must be compatible (!) since 
    \(\extent x = x\) and \(\delta(x,y) = x\wedge y\) and thus
    \(\delta(x,y) = \extent x\wedge\extent y\). Moreover, if 
    \(x_A\) glues \(A\), then 
    \[
        \extent x_A = 
        x_A = 
        \bigvee_{a\in A}\extent a = 
        \bigvee_{a\in A} a
    \]
    and hence \(x_A\) is \(\sup A\). It's easy to verify that the
    supremum is -- if it exists in \(S\) -- the gluing of \(A\). 

    So, the set \(S=\{\bot,a,\top\}\) is gluing complete, 
    obviously. But the singleton \(\sigma(x) = \neg a\wedge x\) 
    (which is the restriction of the singleton represented by 
    \(\neg a\)) has no representing element. 
    \[
        \sigma(\bot) = \bot 
        \qquad
        \sigma(a) = \bot 
        \qquad
        \sigma(\top) = \neg a
    \]
    And so \(S\) is gluing-complete but not Scott-complete.
\end{proof}

\begin{definition}[The \(\quantale Q\)-set of Singletons]
    It is, in general, not possible to form a \(\quantale Q\)-set 
    of singletons over a \(\quantale Q\)-set \(X\). There is some
    trouble in defining \(\delta\) for singletons in a way akin to
    \cite[thm.~2.9.5]{borceux1994handbook3}. This is due to the 
    fact that \(\bigvee_{x}\sigma(x)\tensor\sigma(x)\) not being
    “extent-like” necessarily. 

    We must -- seemingly -- add the axiom of Strength to our 
    quantale \(\quantale Q\) to guarantee we can form a 
    \(\quantale Q\)-set of singletons.
    This enables us to define the following:

    \[
        |\Singletons(X)| = \set{\sigma}{
            \text{\(\sigma\) is a singleton over \(X\)}
        }
    \]
    \[
        \delta(\sigma,\xi) = 
        \bigvee_{x\in X}
            \sigma(x)\tensor\xi(x)
    \]
\end{definition}

\begin{prop}
    Obviously, we must prove this is indeed a \(\quantale Q\)-set.
\end{prop}
\begin{proof}
    Firstly, we know that \(\bigvee_x\sigma(x)\) is idempotent, 
    and hence 
    \begin{align*}
        \bigvee_x\sigma(x) &\leq \bigvee_x\sigma(x)
        \intertext{since \(\quantale Q\) is hypothesized to be 
        strong, we have}
        \bigvee_x\sigma(x) &\leq\bigvee_x
            \sigma(x)\tensor\sigma(x)
        \intertext{and since \(\quantale Q\) is semi-cartesian}
        \tag{semi-cartesian}
        \bigvee_x\sigma(x)\tensor\sigma(x)
            &\leq\bigvee_x\sigma(x)
    \end{align*}
    and therefore, \(\extent \sigma = \bigvee_x\sigma(x)\). Now we
    move onto the \(\quantale Q\)-set axioms. \(\delta\) as defined 
    is obviously symmetric; let's show extensionality:
    \begin{align*}
        \delta(\sigma,\xi)\tensor\extent\sigma &= 
        \bigvee_x \sigma(x)\tensor\xi(x)\tensor\extent \sigma
        \intertext{now recall \(\sigma\)'s strictness 
        condition}
        &= \bigvee_x \sigma(x)\tensor\xi(x)\\
        \intertext{Now we need only to show triangle 
        inequality/transitivity}
        \delta(\sigma, \xi)\tensor\delta(\xi,\psi)
        &=	\bigvee_x\bigvee_y
            \sigma(x)\tensor\xi(x)\tensor
            \xi(y)\tensor\psi(y)
        \\&\leq	\bigvee_x\bigvee_y
            \sigma(x)\tensor\delta(x,y)\tensor\psi(y)
        \\&\leq	\bigvee_y
            \sigma(y)\tensor\psi(y)
        \\&=	\delta(\sigma,\psi)
    \end{align*}
\end{proof}

\begin{lemma}
    A relational morphism \(\phi:X\to Y\) induces a functional map 
    \(\check\phi:X\to\Singletons(Y)\) given by
    \[
        \check\phi(x) = \phi(x,\blank)
    \]
\end{lemma}
\begin{proof}
    We ought to show that the map does indeed take every \(x\) to a
    singleton \(\check\phi(x)\) and that the map is indeed a 
    morphism. 
    \begin{align*}
        \delta(x,x')\tensor\check\phi(x)
        &=	\delta(x,x')\tensor\phi(x,\blank)
        \\&\leq	\phi(x',\blank)
        \\\\
        \extent x\tensor\check\phi(x)
        &=	\extent x\tensor\phi(x,\blank)
        \\&=	\phi(x,\blank)
        \\\\
        [\check\phi(x)](y)\tensor[\check\phi(x)](y')
        &=	\phi(x,y)\tensor\phi(x,y')
        \\&\leq	\delta(y,y')
        \\\\
        [\check\phi(x)](y)\tensor
        \bigvee_{y'}[\check\phi(x)](y') 
        &=	\phi(x,y)\tensor\bigvee_y\phi(x,y') 
        \\&=	\phi(x,y)\tensor\extent x
        \\&=	\phi(x,y)
        \\&=	\check\phi(x')
        \intertext{
            Establishing that \(\check\phi(x)\) is always
            a singleton over \(Y\). Now, with regards to it
            being a morphism, just above we have also shown
            that \(\extent\check\phi(x) = \extent x\), so 
            only one axiom remains:
        }
        \delta(\check\phi(x),\check\phi(x'))
        &=	\bigvee_{y}\phi(x,y)\tensor\phi(x',y)
        \\&\geq	\bigvee_{y}
                \phi(x,y)\tensor
                \phi(x,y)\tensor
                \delta(x,x')
        \\&\geq	\delta(x,x')\tensor\bigvee_{y}
                \phi(x,y)\tensor
                \phi(x,y)
        \\&=	\delta(x,x')\tensor\extent x
        \\&=	\delta(x,x')
    \end{align*}
\end{proof}

\begin{remark}
    It is not \emph{necessary} for \(\quantale Q\) to be strong for
    the question of Scott completeness to make sense; but it is not 
    apparent how one may perform Scott-\emph{completion} without 
    that added strength. We will talk about completion now.
\end{remark}

\begin{definition}
    We denote by \(\ScottQSet\) the full subcategory of \(\QSet_f\)
    whose objects are Scott-complete.
\end{definition}

\begin{definition}[Separability]
    A \(\quantale Q\)-set is said to be \emph{separable} when the 
    mapping \(x\mapsto\delta(\blank, x)\) is \emph{injective}; 
    \lat{ie.} it is the completeness condition sans surjectivity.
\end{definition}

\begin{prop}
    Separability is equivalent to Extensionality.
\end{prop}
\begin{proof}
    Immediate from what was proven in the definition of 
    \(\delta\)-equivalence.
\end{proof}

\subsection{Scott Completion}
Scott-completion is, unsurprisingly, going to be a reflector functor of
the from \(\QSet\) to the full subcategory \(\ScottQSet\); in this 
subsection we therefore will show that \(\Singletons(X)\) is indeed 
Scott-complete and that object assignment can be made functorial -- in 
such a way to get a left adjoint to the fully faithful subcategory 
inclusion.

\begin{lemma}[\cite{borceux1994handbook3} cf. p. 158]
    If \(\sigma_x\) denotes \(\delta(\blank,x)\), then 
    \[
        \delta(\sigma_x,\xi) = \xi(x)
    \]
    a result akin to the Yoneda lemma, wherein 
    \[
        \hom(\YonedaEmbedding(x), F)\iso F(x)
    \]
\end{lemma}
\begin{proof}
    \begin{align*}
        \delta(\sigma_x,\xi) 
        &=	\bigvee_{y}\delta(y,x)\tensor\xi(y)
        \\&\leq	\xi(x) 
        \\&=	\extent x\tensor \xi(x)
        \\&=	\delta(x, x) \tensor \xi(x)
        \\&\leq	\bigvee_{z}\delta(z, x) \tensor \xi(z)
        \\&=	\delta(\sigma_x,\xi)
    \end{align*}
\end{proof}
\begin{remark}
    Notice that this means 
    \(\delta(\sigma_x,\sigma_y) = \delta(x,y)\)
    and thus the mapping \(x\mapsto\sigma_x\) preserves 
    \(\delta\) and hence is a morphism -- it is in fact a regular
    monomorphism.
\end{remark}

\begin{prop}
    \(\Singletons(X)\) is, indeed, Scott-complete.
\end{prop}
\begin{proof}
    In this proof, we shall refer to members of \(\Singletons(X)\)
    as \(1\)-singletons and use lowercase \(\xi,\psi,\phi\)
    and refer to members of \(\Singletons\circ\Singletons(X)\) as 
    \(2\)-singletons and use uppercase \(\Sigma\) \lat{etc.}

    To go about doing this we must show two things: that two 
    \(1\)-singletons cannot represent the same \(2\)-singleton; and
    that every \(2\)-singleton comes from a \(1\)-singleton.

    For the purpose of showing that \(\xi\mapsto\sigma_\xi\) is 
    injective, suppose that \(\xi\) and \(\psi\) are such that 
    \(\sigma_\xi=\sigma_\psi\). Let \(x\in X\) and realize
    \[
        \xi(x) = \delta(\sigma_x, \xi) = \sigma_\xi(\sigma_x)
        = \sigma_\psi(\sigma_x) = \delta(\sigma_x,\psi) = 
        \psi(x)
    \]
    thence they are extensionally equal.

    Now take a \(2\)-singleton \(\Sigma\), we shall define a 
    \(1\)-singleton that we hope represents it; namely:
    \[
        \phi(x) = \Sigma(\sigma_x)
    \]
    one way to become convinced that \(\phi\) is indeed a singleton
    over \(X\) is that \(\sigma_{\blank}\) preserves \(\delta\) and 
    thus \(\Sigma\circ\sigma_{\blank}\) is “more or less” the 
    restriction of \(\Sigma\) to the set of representable 
    singletons. The one axiom that doesn't obviously hold is 
    strictness -- which we show after the following deduction:
    \begin{align*}
        \Sigma(\xi) 
        &=	\extent\xi\tensor\Sigma(\xi)
        \\&=	\bigvee_x\xi(x)\tensor\xi(x)\tensor\Sigma(\xi)
        \\&=	\bigvee_x
            \xi(x)\tensor
            \delta(\sigma_x,\xi)\tensor
            \Sigma(\xi)
        \\&\leq	\bigvee_x
            \xi(x)\tensor
            \Sigma(\sigma_x)
        \\&=	\delta(\xi,\phi)
    \end{align*}

    Hence, if we can show strictness, we are done. So take 
    \(x\in X\) which we ought to show enjoys the property that
    \[
        \phi(x)\tensor\bigvee_y\phi(y) = \phi(x)
    \]
    it is clear that it would suffice if we proved that 
    \(\extent\Sigma \leq \extent \phi\). For that would suffice for
    them to be indeed equal; and them being equal mean we could 
    then argue:
    \[
        \phi(x)\tensor\bigvee_y\phi(y) = 
        \phi(x)\tensor\extent \phi = 
        \Sigma(\sigma_x)\tensor\extent\Sigma = 
        \Sigma(\sigma_x) = \phi(x)
    \]
    since \(\Sigma\) was taken to be strict. And thus we shall 
    prove instead simply that \(\extent\Sigma\leq\extent\phi\).

    \begin{align*}
        \extent\Sigma 
        &=	\bigvee_\xi\Sigma(\xi)\tensor\Sigma(\xi)
        \\&=	\bigvee_\xi\Sigma(\xi)
        \tag{Strength}
        \\&=	\bigvee_\xi\Sigma(\xi)\tensor\extent\xi
        \\&=	\bigvee_x\bigvee_\xi
                \Sigma(\xi)\tensor\xi(x)
        \\&=	\bigvee_x\bigvee_\xi
                \Sigma(\xi)\tensor
                \delta(\xi,\sigma_x)
        \\&\leq	\bigvee_x\Sigma(\sigma_x)
        \\&=	\bigvee_x
                \Sigma(\sigma_x)\tensor\Sigma(\sigma_x)
        \tag{Strength}
        \\&=	\extent\phi
    \end{align*}
\end{proof}

Having then proved that \(\Singletons(X)\) is indeed Scott-complete, 
we now must move into the associated completion functor.

\begin{definition}
    We extend the action of \(\Singletons\) by acting on functional
    morphisms as follows:
    \begin{align*}
        f : X&\to Y\\
        \\
        \Singletons f : \Singletons(X)&\to\Singletons(Y)\\
        \xi&\mapsto
            \bigvee_{x\in X}
            \delta(f(x),\blank)\tensor\xi(x)
    \end{align*}
    Morally speaking, if \(\xi(x)\) measures how much \(x\) is akin
    to some abstract point, then \((\Singletons f)(\xi)\) measures 
    how close a point in \(y\) is to the “image” of that abstract 
    point under \(f\).
\end{definition}

\begin{prop}
    Associated to the previous definition, we must prove that 
    \(\Singletons f(\xi)\) is indeed a singleton; that this 
    assignment is a morphism; and that the action on morphisms is 
    functorial.
\end{prop}
\begin{proof}
    First of all, for notational sanity denote 
    \(\Singletons f(\xi)\) by \(\xi_f\). 
    \begin{align*}
        \xi_f(y)\tensor\extent y 
        &=	\bigvee_x\delta(f(x),y)\tensor\xi(x)
            \tensor\extent y 
        \\&= 	\bigvee_x\delta(f(x),y)\tensor\xi(x)
        \\&=	\xi_f(y)\\
        \\
        \xi_f(y)\tensor\delta(y,y')
        &=	\bigvee_x\delta(f(x),y)\tensor
            \delta(y,y')\tensor\xi(y)
        \\&\leq	\bigvee_x\delta(f(x),y')\tensor
            \tensor\xi(y)
        \\&=	\xi_f(y')\\
        \\
        \xi_f(y)\tensor\xi_f(y') 
        &=	\bigvee_{x,x'}
            \delta(f(x), y)\tensor
            \delta(f(x'), y')\tensor
            \xi(x)\tensor\xi(x')
        \\&\leq	\bigvee_{x,x'}
            \delta(f(x), y)\tensor
            \delta(f(x'), y')\tensor
            \delta(x,x')
        \\&\leq	\bigvee_{x,x'}
            \delta(f(x), y)\tensor
            \delta(f(x'), y')\tensor
            \delta(f(x),f(x'))
        \\&\leq	\bigvee_{x'}
            \delta(f(x'), y)\tensor
            \delta(f(x'), y')
        \\&\leq	\delta(y,y')\\
        \\
        \xi_f(y)\tensor\bigvee_{y'}\xi_f(y')
        &=	\bigvee_{y'}\bigvee_{x'}
                \xi_f(y)\tensor
                \delta(f(x'),y)\tensor
                \xi(x')
        \\&=	\bigvee_{x'}
                \xi_f(y)\tensor
                \delta(f(x'),f(x'))\tensor
                \xi(x')
        \\&=	\bigvee_{x'}
                \xi_f(y)\tensor
                \extent x'\tensor
                \xi(x')
        \\&=	\bigvee_{x'}
                \xi_f(y)\tensor
                \xi(x')
        \\&=	\bigvee_{x}
                \delta(f(x), y)\tensor
                \xi(x)\tensor
                \bigvee_{x'}
                    \xi(x')
        \\&=	\bigvee_{x}
                \delta(f(x), y)\tensor
                \xi(x)
        \\&=	\xi_f(y)
    \end{align*}

    The second thing we ought to do is showing that the 
    map \(\Singletons f\) is a morphism; that is rather 
    easy actually. Let \(\psi_f\) denote the image of \(\psi\)
    under \(\Singletons f\) as we had done with \(\xi\).
    \begin{align*}
        \delta(\xi_f,\psi_f) 
        &=	\bigvee_y\xi_f(y)\tensor\psi_f(y)
        \\&=	\bigvee_y
            \bigvee_{x, x'}
                \delta(f(x),y)
                \tensor
                \delta(f(x'),y)
                \tensor
                \xi(x)
                \tensor
                \phi(x')
        \\&\geq	\bigvee_y
            \bigvee_{x}
                \delta(f(x),y)
                \tensor
                \delta(f(x),y)
                \tensor
                \xi(x)
                \tensor
                \phi(x)
        \tag{by letting \(x=x'\)}
        \\&\geq	\bigvee_{x}
                \extent f(x)
                \tensor
                \extent f(x)
                \tensor
                \xi(x)
                \tensor
                \phi(x)
        \tag{by letting \(y=f(x)\)}
        \\&=	\bigvee_{x}
                \xi(x)
                \tensor
                \phi(x)
        \\&=	\delta(\xi, \phi)\\
        \\
        \xi_f(y)
        \tag{semi-cartesian}
        &\geq	\xi_f(y)\tensor\bigvee_{y'}\xi_f(y')
        \\&=	\bigvee_{y'}\bigvee_{x'}
                \xi_f(y)\tensor
                \delta(f(x'), y')\tensor
                \xi(x')
        \\&\geq	\bigvee_{x'}
                \xi_f(y)\tensor
                \delta(f(x'),f(x'))\tensor
                \xi(x')
        \tag{by letting \(y'=f(x')\)}
        \\&=	\bigvee_{x'}
                \xi_f(y)\tensor
                \extent f(x')\tensor
                \xi(x')
        \\&=	\xi_f(y)\tensor
            \bigvee_{x'}
                \xi(x')
        \\&=	\bigvee_{x}
                \delta(f(x),y)\tensor\xi(x)\tensor
            \bigvee_{x'}
                \xi(x')
        \\&=	\bigvee_{x}
                \delta(f(x),y)\tensor\xi(x)
        \\&=	\xi_f(y)
    \end{align*}

    The final problem to be unravelled is functoriality. We must 
    show that \(\id\) is sent to \(\id\) and composition is 
    preserved. \(\Singletons\id\) is defined as
    \begin{align*}
        [\Singletons\id(\xi)](x)
        &=	\bigvee_{x'}\delta(x,x')\tensor\xi(x')
        \\&\geq	\extent x\tensor\xi x = \xi(x)
        \\&\geq	\delta(x', x)\tensor\xi(x), \ \forall x'
    \end{align*}
    and hence by universality of the supremum, we have the 
    desired equality: 
    \[
        (\Singletons\id)(\xi) = \xi
    \]
    Now, take \(X\xrightarrow f Y\xrightarrow g Z\)
    \begin{align*}
        [\Singletons(f\circ g)(\xi)](z)
        &=	\bigvee_x\delta(g\circ f(x), z)\tensor\xi(x)
        \\&=	\bigvee_x
            \bigvee_y
                \delta(f(x), y)\tensor\xi(x)
                \tensor
                \delta(g(y), z)
        \\&=	\bigvee_y
            \delta(g(y),z)\tensor
            \left(
                \bigvee_x\delta(f(x),y)\tensor\xi(x)
            \right)
        \\&=	\bigvee_y
            \delta(g(y),z)\tensor
            [\Singletons f(\xi)](y)
        \\&=	[(\Singletons g)((\Singletons f)(\xi))](z)
        \\&=	[(\Singletons g)\circ(\Singletons f)](z)
    \end{align*}
    Thus, we have shown that \(\Singletons\) is a functor from
    \(\QSet_f\) to \(\ScottQSet\).
\end{proof}

\begin{theorem}
    We aim to show that \(\Singletons\) indeed a completion, in the 
    sense that it is the left adjoint to an full inclusion functor 
    -- of \(\ScottQSet\); yielding that if \(K\) is Scott-complete,
    \[
        \QSet(X, K)\iso\QSet(\Singletons X, K)
    \]
\end{theorem}
\begin{proof}
    Again, we shall proceed by finding units and counits. Once more
    let us rename \(\Singletons\) by \(L\) dub the inclusion 
    functor simply \(R\), after their role in the adjunction. We 
    ought to find \(\eta\) and \(\epsilon\) as we had in 
    Gluing-completion proof.

    The unit, \(\eta\) is straightforward, it is the associated 
    singleton morphism \(x\mapsto\delta(\blank, x)\) -- which we 
    know to be a morphism already, and hence we need only to show 
    naturality from \(\id\) to \(R\circ L\).
    \[\begin{tikzcd}
        X &&& {\Singletons(X)} \\
        & x & {\delta(\blank, x)} \\
        & {f(x)} & {?} \\
        Y &&& {\Singletons(Y)}
        \arrow["{\Singletons(f)}"{description}, 
            from=1-4, to=4-4]
        \arrow["\eta", from=4-1, to=4-4]
        \arrow["\eta", from=1-1, to=1-4]
        \arrow["f"{description}, from=1-1, to=4-1]
        \arrow[maps to, from=2-2, to=3-2]
        \arrow[maps to, from=3-2, to=3-3]
        \arrow[maps to, from=2-2, to=2-3]
        \arrow[maps to, from=2-3, to=3-3]
    \end{tikzcd}\]
    \[
        \Singletons f(\delta(\blank,x)) = 
        \bigvee_{x'}\delta(\blank,f(x'))\tensor\delta(x', x)=
        \delta(\blank, f(x))
    \]
    and hence we have naturality. It should be remarked that 
    \(\Singletons f(\sigma_x) = \sigma_{f(x)}\) as just proved 
    above.

    The counit is also rather straightforward; we start with a 
    Scott complete \(X\), we forget all about completeness and 
    perform a completion. Since we know that \(x\mapsto\sigma_x\)
    is bijective (by definition of completeness) and that this 
    mapping preserves \(\delta\), we know that \(X\) is isomorphic
    to \(\Singletons(X)\). So we need only show that this map
    is a \emph{natural} isomorphism. Here \(\epsilon\) takes a 
    singleton (which is always representable) to its unique 
    representing element:
    \[\begin{tikzcd}
        {\Singletons(X)} &&& X \\
        & {\delta(\blank, x)} & x \\
        & {\delta(\blank, f(x))} & {?} \\
        {\Singletons(Y)} &&& Y
        \arrow["{\Singletons(f)}"{description}, from=1-1, 
            to=4-1]
        \arrow["\epsilon", from=4-1, to=4-4]
        \arrow["\epsilon", from=1-1, to=1-4]
        \arrow["f"{description}, from=1-4, to=4-4]
        \arrow[maps to, from=2-2, to=3-2]
        \arrow[maps to, from=3-2, to=3-3]
        \arrow[maps to, from=2-2, to=2-3]
        \arrow[maps to, from=2-3, to=3-3]
    \end{tikzcd}\]
    Which is obviously commutative. Recapping the definitions,

    \begin{align*}
        \eta_X: X&\to\Singletons(X)\\
        x&\mapsto \delta(\blank, x)\\
        \\
        L\eta : L
        &\xrightarrow{\phantom{(L\eta)_X}}
        L\circ R\circ L\\
            (X\in\QSet)\mapsto( 
                \xi
                &\xmapsto{(L\eta)_X}
                \bigvee_x
                    \delta(\blank, \eta(x))\tensor
                    \xi(x)
            )
        \intertext{As we know from the Yoneda-esque lemma, the 
        above can be simplified into the following:}
        L\eta : L
        &\xrightarrow{\phantom{(L\eta)_X}}
        L\circ R\circ L\\
            (X\in\QSet)\mapsto( 
                \xi
                &\xmapsto{(L\eta)_X}
                \bigvee_x
                    \blank(x)\tensor\xi(x)
            )
        \intertext{
            which we know to be simply 
            \(\delta(\blank,\xi)\). As for \(\epsilon L\), 
            we first recall what is the representing 
            element of a \(2\)-singleton \(\Sigma\). It is 
            what we had called \(\phi\) -- which was the 
            action of \(\Sigma\) on representable 
            singletons, or more helpfully, the composite 
            \(\Sigma(\sigma_{\blank})\) or even better: 
            \(\Sigma\circ\eta\).
        }
        \\
        \epsilon L: L\circ R\circ L
        &\xrightarrow{\phantom{(\epsilon R)_X}}
        L\\
            (X\in\QSet)\mapsto(
                \Sigma
                &\xmapsto{(\epsilon L)_X}
                \Sigma\circ\eta
            )\\
        \\
        R\epsilon: R\circ L\circ R
        &\xrightarrow{\phantom{R\epsilon)_X\ }}
        R\\
            (X\in\GluingQSet)\mapsto(
                \delta(\blank, x)
                &\xmapsto{(R\epsilon)_X}
                x
            )\\
        \\
        \eta R: R
        &\xrightarrow{\phantom{\eta R)_X\ }}
        R\circ L\circ R\\
            (X\in\GluingQSet)\mapsto(
                x
                &\xmapsto{(\eta R)_X}
                \delta(\blank, x)
            )
    \end{align*}
    Tracing \(x\) along \((R\epsilon)\circ(\eta R)\) we find
    \[
        x
        \xmapsto{(\eta R)_X}
        \delta(\blank, x)
        \xmapsto{(R\epsilon)_X)}
        x
    \]
    and now, tracing \(\xi\) along \((\epsilon L)\circ(L\eta)\),
    we get
    \[
        \xi
        \xmapsto{(L\eta)_X}
        \delta(\blank,\xi)
        \xmapsto{(\epsilon L)_X}
        \delta(\sigma_{\blank},\xi)
    \]
    which, again, by the Yoneda-esque lemma is \(\xi(\blank)\) 
    which is extensionally equal to \(\xi\), of course. And thus,
    our two composites were equal to the appropriate identities;
    thus establishing the adjunction. 

    We already knew the counit was an isomorphism, but since it 
    is the couint of an adjunction with a fully faithful right 
    adjoint that would have proved it as well.
\end{proof}
	\section{Connection between Completeness Conditions}


\begin{theorem} 
    Let \((X, \delta) \) be an extensional \(\quantale Q\)-set. 
    The following conditions are equivalent:
    \begin{enumerate}
        \item \((X, \delta) \) is Scott-complete.
        \item  \((X, \delta) \) is gluing-complete and, for each singleton \(\sigma\)  over \((X, \delta) \), it holds the condition below:\\
        \((*_\sigma)\)  \ For each \(x \in X\), there is \(y \in X\) such that \(\sigma(x) \leq \sigma(y) = \extent y\).
    \end{enumerate}
\end{theorem}

We have already mentioned that Scott-completeness implies 
Gluing-completeness; and hence there is another fully faithful 
subcategory inclusion at play, forming the triangle:
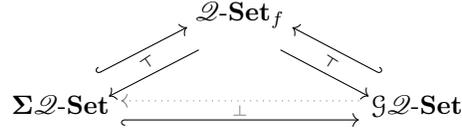
\[\begin{tikzcd}
    & {\QSet_f} \\
    \ScottQSet && \GluingQSet
    \arrow[""{name=0, anchor=center, inner sep=0}, shift left=2, 
        hook, from=2-1, to=1-2]
    \arrow[""{name=1, anchor=center, inner sep=0}, shift right=2, 
        from=1-2, to=2-3]
    \arrow[""{name=2, anchor=center, inner sep=0}, shift right=2, 
        hook', from=2-1, to=2-3]
    \arrow[""{name=3, anchor=center, inner sep=0}, shift left=2, 
        from=1-2, to=2-1]
    \arrow[""{name=4, anchor=center, inner sep=0}, shift right=2, 
        hook', from=2-3, to=1-2]
    \arrow[""{name=5, anchor=center, inner sep=0}, shift right=1, 
        color={rgb,255:red,128;green,128;blue,128}, dotted, 
        from=2-3, to=2-1]
    \arrow["\dashv"{anchor=center, rotate=124}, draw=none, from=3, 
        to=0]
    \arrow["\dashv"{anchor=center, rotate=58}, draw=none, from=1, 
        to=4]
    \arrow["\dashv"{anchor=center, rotate=-90}, 
        color={rgb,255:red,128;green,128;blue,128}, draw=none, 
        from=5, to=2]
\end{tikzcd}\]
The question, visually posed in the diagram by the gray dotted arrow,
is “is there a functor left adjoint to that lonely inclusion?” The,
answer, not surprisingly, is yes. And it is given by the only composite
going in that direction. Now, this isn't very spectacular, and it comes 
directly from \(\hom\) isomorphism and the categories involved being 
full \lat{etc.} 

The more interesting categorical property is the characterization of 
Scott-completeness; which we shall work towards in this next 
subsection.

\subsection{Scott-Completeness and Relational Morphisms}
Scott-completeness is deeply tied to relational morphisms -- which 
might be surprising since we have deliberately not dealt with 
\(\QSet_r\) in the context of completeness thus far. The reason is that
singleton completeness is \emph{too} deeply tied to relational 
morphisms.

\begin{theorem}
    \[\ScottQSet_r\equiv\QSet_r\]
\end{theorem}
\begin{proof}
    We proceed by providing a functor fully faithful functor that 
    is essentially surjective; in this case, it is easier to show
    that the full subcategory inclusion is essentially surjective,
    since by necessity is fully faithful.

    This then amounts to us provinding an isomorphism between any
    object and a Scott-complete one. In this case we have an 
    obvious candidate, that being \(\Singletons(X)\). Define
    \begin{align*}
        \phi:X &\to \Singletons(X)\\
            (x,\xi)&\mapsto\xi(x)\\
        \\
        \phi^{-1}:\Singletons(X)&\to X\\
            (\xi, x)&\mapsto\xi(x)
    \end{align*}
    We ought to prove that both of them are morphisms first, but 
    should the reader grant us a moratorium on that we could 
    argue that
    \begin{align*}
        \phi\circ\phi^{-1}(\xi, \psi) 
        &=	\bigvee_{x}\xi(x)\tensor\psi(x)
        \\&=	\delta(\xi, \psi)\\
        \\
        \phi^{-1}\circ\phi(x,x') 
        &=	\bigvee_{\xi}\xi(x)\tensor\xi(x')
        \\&\geq	\delta(x,x)\tensor\delta(x,x')
        \\&=	\delta(x,x')
    \end{align*}
    hence, their composite is pointwise greater than the identity; 
    but we know that under strong quantales, the pointwise order 
    on \(\hom\)s is in fact discrete. So they are actually equal.

    It remains to be seen that as defined, both \(\phi\) and its 
    supposed inverse are actually morphisms. To it, then.
    It should be obvious why the extensionality axioms should hold:
    \begin{align*}	
        \xi(x)\tensor\extent\xi 
        &=	\xi(x)\tensor\bigvee_{x'}\xi(x') 
        \\&=	\xi(x) 
        \\&= 	\xi(x)\tensor\extent x
        \intertext{The \(\delta\) laws also come for free:}
        \delta(x,x')\tensor\xi(x)
        &\leq\xi(x)
        \tag{subset axiom}
        \\
        \\
        \delta(\xi,\psi)\tensor\xi(x) 
        &=	\delta(\xi,\psi)\tensor\delta(\sigma_x,\xi)
        \\&\leq	\delta(\sigma_x,\psi)
        \\&=	\psi(x)
        \intertext{
            The singleton condition holds trivially as well
            -- due to how \(\delta(\xi,\psi)\) is defined.
        }
        \phi(x,\xi)\otimes\phi(x,\psi) 
        &=	\xi(x)\tensor\psi(x) 
        \\&\leq	\bigvee_{x}\xi(x)\tensor\psi(x)
        \\&=	\delta(\xi,\psi)
        \\\\
        \phi^{-1}(\xi, x)\tensor\phi^{-1}(\xi, x') 
        &=	\xi(x)\tensor\xi(x')
        \tag{\(\xi\) is a singleton}
        \\&\leq	\delta(x,x')
        \intertext{remaining only to show strictness}
        \bigvee_\xi\phi(x,\xi) 
        &=	\bigvee_{\xi}\xi(x)
        \tag{now let \(\xi = \sigma_x\)}
        \\&\geq	\delta(x,x) = \extent x
        \\&\geq	\xi(x)~\forall\xi
        \intertext{
            and hence \(\phi\) is indeed strict. Now for 
            \(\phi^{-1}\):
        }
        \bigvee_{x}\phi^{-1}(\xi,x) 
        &=	\bigvee_x\xi(x)
        \\&=	\extent\xi
    \end{align*}
\end{proof}

\begin{theorem}
    \[
        \ScottQSet_f\iso\ScottQSet_r
    \]

    Recall the underlying graph functor \(\Graph\); we propose that
    its restriction to \(\ScottQSet_f\) has an inverse, and we give
    it explicitly now: 
    \[
        \Functional:\ScottQSet_r\to\ScottQSet_f
    \] 
    is a functor that “remembers” that relational morphisms between
    \(\quantale Q\)-sets are always functional. It's action on 
    objects is the identity. As we had mentioned in the previous 
    subsection, relational morphisms \(\phi:X\to Y\) induce a 
    functional morphism \(\check\phi:X\to\Singletons(Y)\) given by 
    \(x\mapsto\phi(x,\blank)\). 

    Since \(Y\) is Scott complete, it is naturally isomorphic to 
    \(\Singletons(Y)\) and hence each \(\check\phi(x)\) corresponds
    to some \(y_{\check\phi(x)}\) which is the unique point in 
    \(y\) representing \(\check\phi(x)\). We claim that this 
    mapping is a functional morphism \(X\to Y\) (easy) and that 
    this action on morphisms is functorial, and assign it to 
    \(\Functional\); we further claim that \(\Functional\) is the 
    inverse of \(\Graph\restriction_{\ScottQSet_f}^{\ScottQSet_r}\)
    -- which is a nasty thing to typeset, so we'll just say 
    \(\Graph\).
\end{theorem}
\begin{proof}
    There is quite a lot to unpack; formally, we say that if we 
    define \(\Functional(\phi) = \epsilon\circ\check\phi\), then
    this action is both functorial and makes \(\Functional\) and 
    \(\Graph\) mutual inverses --- where \(\epsilon\) is the counit
    of the adjunction, which takes singletons over Scott-complete 
    \(\quantale Q\)-sets to their representing elements.

    Let's write, for a given \(\phi\), \(\Functional(\phi)\) as 
    \(f_\phi\) for simplicity sake. First, note that 
    \(\phi(x,f_\phi(x))=\extent x\) -- as 
    \(\phi(x,\blank) = \delta(\blank, f_\phi(x))\) and 
    \(\extent f_\phi(x) = \extent x\). Obviously, \(\Functional\) 
    preserves identities, as
    \begin{align*}
        f_{\id}
        &=	\epsilon(x\mapsto(\blank,x)) 
        \\&=	x\mapsto x 
        \\&=	\id
        \intertext{
            functoriality can now be seen to hold: given 
            \(x\), consider the singleton
        } 
        \bigvee_y\phi(x,y)\tensor\psi(y,\blank)
        &\geq	\phi(x,f_\phi(x))\tensor\psi(f_\phi(x),\blank)
        \\&=	\extent f_\phi(x)\tensor\psi(f_\phi(x),\blank)
        \\&=	\psi(f_\phi(x),\blank)
        \\&\geq	\psi(y,\blank)\tensor\delta(f_\phi(x),y)
        \tag{for any \(y\)}
        \\&=	\psi(y,\blank)\tensor\phi(x,y)
        \intertext{
            and hence 
            \([\psi\circ\phi](x,\blank) = 
            \psi(f_\phi(x),\blank)\), from which we can 
            conclude that
        }
        f_{\psi\circ\phi} &= f_{\psi}\circ f_\phi
    \end{align*}

    We ought to show that \(\Graph\circ\Functional(\phi) = \phi\) 
    and \(\Functional\circ\Graph(f) = f\). Let's go about doing it 
    in the order we've written. Take some relational morphism 
    \(\phi:X\to Y\) and let's consider the following:
    \begin{align*}
        [\Graph\circ\Functional(\phi)](x,y) 
        &=	[\Graph(\epsilon\circ\check\phi)](x,y)
        \\&=	\delta(\epsilon\circ\check\phi(x), y)
        \\&=	\delta(\check\phi(x), \sigma_y)
        \\&=	\check\phi(x)(y)
        \\&=	\phi(x,y)
        \intertext{
            Now take some functional morphism \(f:X\to Y\)
        }
        [\Functional\circ\Graph(f)](x) 
        &=	[\Functional(\delta(f(
                \blank[1]), 
                \blank[2]
            ))](x)
        \\&=	[f(\blank[1])](x) 
        \\&=	f(x)
    \end{align*}
    And therefore, we have shown that the composites are indeed 
    actually the identities and \(\Functional\) witnesses the fact
    that \(\ScottQSet_r\iso\ScottQSet_f\).
\end{proof}

\begin{lemma}
    \[\Graph\circ\Singletons\iso\Graph\]
    Where \(\Graph\circ\Singletons\) really ought to be the composite
    of \(\Singletons\) with the appropriate restriction of \(\Graph\) 
    to the full subcategory that is the image of \(\Singletons\). Or,
    alternatively, \(\Graph\circ i \circ\Singletons\).
\end{lemma}
\begin{proof}
    We must provide an \(\alpha_X\) for every \(X\) such that
    \[\begin{tikzcd}
        X & {i\circ\Singletons(X)} \\
        Y & {i\circ\Singletons(X)}
        \arrow["\alpha", from=1-1, to=1-2]
        \arrow["\alpha", from=2-1, to=2-2]
        \arrow["{\Graph(f)}"{description}, from=1-1, to=2-1]
        \arrow["{\Graph\circ i\circ\Singletons(f)}"{description}, 
            from=1-2, to=2-2]
    \end{tikzcd}\]
    which simplifies to finding \(\alpha_X\) such that the following
    commutes, thus justifying our lack of precision in the “clean” 
    statement of the lemma:
    \[\begin{tikzcd}
        X & {\Singletons(X)} \\
        Y & {\Singletons(X)}
        \arrow["\alpha", from=1-1, to=1-2]
        \arrow["\alpha", from=2-1, to=2-2]
        \arrow["{\Graph(f)}"{description}, from=1-1, to=2-1]
        \arrow["{\Graph\circ\Singletons(f)}"{description}, 
            from=1-2, to=2-2]
    \end{tikzcd}\]
    if we then let \(\alpha\) be the relational isomorphism we 
    already know exists between \(X\) and \(\Singletons(X)\) there 
    is a good chance naturality holds by magic. We need to verify 
    that 
    \[
        \alpha\circ[\Graph(f)] = 
        [\Graph\circ\Singletons(f)]\circ
        \alpha
    \]
    We do so by \emph{cheating}, and proving instead that 
    \[
        \alpha\circ[\Graph(f)] \geq 
        [\Graph\circ\Singletons(f)]\circ
        \alpha
    \]
    and then we remember that for relational morphisms \(\leq\) 
    \emph{is} \(=\) because \(\quantale Q\) is \textbf{strong}.

    \begin{align*}
        [\Graph\circ\Singletons(f)]\circ\alpha(x,\psi)
        &=	\bigvee_\xi
                \xi(x)\tensor
                \delta(\Singletons f(\xi),\psi)
        \\&=	\bigvee_\xi
                \xi(x)\tensor
            \bigvee_{y}
            \bigvee_{x}
                \delta(f(x),y)
                \tensor
                \xi(x)
                \tensor
                \psi(y)
        \\&=	\bigvee_{y}
            \bigvee_{x'}
            \bigvee_\xi
                \xi(x)\tensor
                \delta(f(x'),y)
                \tensor
                \xi(x')
                \tensor
                \psi(y)
        \\&\leq	\bigvee_y
            \bigvee_{x'}
                \delta(x,x')
                \tensor
                \delta(f(x'),y)
                \tensor
                \psi(y)
        \\&\leq	\bigvee_y
            \bigvee_{x'}
                \delta(f(x),f(x'))
                \tensor
                \delta(f(x'),y)
                \tensor
                \psi(y)
        \\&\leq	\bigvee_y
                \delta(f(x),y)
                \tensor
                \psi(y)
        \\&=	\alpha\circ[\Graph(f)](x,\psi)
    \end{align*}
\end{proof}

Now we are ready to state a more categorical description of 
Scott-completeness, merely in terms of representable functors 
\lat{etc.} 
\begin{theorem}
    Let \(\restriction\Graph\) denote the funtor between the 
    presheaf categories of \(\QSet_r\) and \(\QSet_f\) as below:
    \[\begin{tikzcd}
        {\QSet_f} & {\QSet_r} \\
        {\PSh(\QSet_f)} & {\PSh(\QSet_r)}
        \arrow[""{name=0, anchor=center, inner sep=0}, 
            "\Graph", from=1-1, to=1-2]
        \arrow["\PSh"{description}, maps to, from=1-1, to=2-1]
        \arrow["\PSh"{description}, maps to, from=1-2, to=2-2]
        \arrow[""{name=1, anchor=center, inner sep=0}, 
            "\restriction\Graph", from=2-2, to=2-1]
        \arrow["\PSh"{description}, shorten <=4pt, 
            shorten >=4pt, maps to, from=0, to=1]
    \end{tikzcd}\]
    Given by precomposing presheaves with the functor \(\Graph\)
    so as to change their domains appropriately. Then 
    \[
        \left[
            \YonedaEmbedding_f(X)\iso
            \YonedaEmbedding_r(X)\restriction\Graph
        \right]
        \iff
        X ~ \text{is Scott-complete}
    \]
\end{theorem}
\begin{proof}
    If \(X\) is indeed complete, then we know that the following 
    are natural isomorphisms in \(Y\)
    \begin{align*}
        \YonedaEmbedding_f(X)
        &=	\YonedaEmbedding_f(R(X))
        \\&\iso	\YonedaEmbedding_f(X)\circ\Singletons
        \tag{adjunction}
        \intertext{and, as both entries are now complete,}
        &=	\ScottQSet_f(\Singletons(\blank),X)
        \intertext{
            We also know that \(\Graph\) is an isomorphism 
            between the subcategories \(\ScottQSet_f\) and 
            \(\ScottQSet_r\); in particular we have
            \[
                \ScottQSet_r(
                    \Graph(\blank[1]),
                    \Graph(\blank[2])
                )\iso
                \ScottQSet_f(
                    \blank[1],
                    \blank[2]
                )
            \]
            if we fix \(\blank[2] = X\), we can forget 
            about the \(\Graph\) on the second coordinate 
            (because it's there just to act on morphisms),
            thus getting
            \[
                \ScottQSet_r(\Graph(\blank), X)\iso
                \ScottQSet_f(\blank, X)
            \]
            Resuming our isomorphism chain:
        }
        \ScottQSet_f(\Singletons(\blank),X)
        &\iso	\ScottQSet_r(\Graph\circ\Singletons(\blank), X)
        \\&=	\YonedaEmbedding_r(X)\circ\Graph\circ\Singletons
        \intertext{
            For our claim to hold, we would need a natural
            isomorphism \(\Graph\to\Graph\circ\Singletons\)
            -- the reader should be happy to remember that 
            this isomorphism has already been proven in an
            earlier lemma. We can then append one last 
            isomorphism then:
        }
        \YonedaEmbedding_r(X)\circ\Graph\circ\Singletons
        &\iso	\YonedaEmbedding_r\circ\Graph
    \end{align*}
    We can then put it all together and show there is some 
    twisted little natural isomorphism 
    \[
        \YonedaEmbedding_f(X)\iso
        \YonedaEmbedding_r\circ\Graph = 
        \YonedaEmbedding_r(X)\restriction\Graph
    \]

    The converse implication is, thankfully, easier; suppose there
    is indeed an isomorphism as above; we ought to show that \(X\)
    was Scott-complete to begin with. 
    \begin{align*}	
        \YonedaEmbedding_f(i\circ\Singletons(X))
        &\iso	\YonedaEmbedding_r(i\circ\Singletons(X))\circ\Graph
        \tag{just proven}
        \\&\iso	\YonedaEmbedding_r(X)\circ\Graph
        \tag{\(X\iso_r i\circ\Singletons(X)\)}
        \\&\iso	\YonedaEmbedding_f(X)
        \tag{hypothesis}
        \intertext{thence, by Yoneda we obtain that}
        \Singletons(X)&\iso_f X
    \end{align*}
    But Scott-completeness is invariant under functional 
    isomorphisms. So \(X\) was Scott-complete to begin with and
    thus ends our proof.
\end{proof}

We then give a useful tool to show Scott-completeness in terms
of an object's representable functor. Namely
\begin{theorem}
    \[
        X \ \text{is Scott-complete}
        \iff
        \QSet_f(\blank, X)\iso\QSet_f(\Singletons(\blank), X)
    \]
    here again we actually wrote something slightly informal
    and we really mean that, \(\QSet_f(i\circ\Singletons(\blank), X)\)
    is isomorphic to \(\QSet_f(\blank, X)\). We can rewrite it as
    \[
        \YonedaEmbedding_f(X) \iso 
        \YonedaEmbedding_f(X) \circ i \circ \Singletons
    \]
\end{theorem}
\begin{proof}
    Scott-completeness implies the existence of the isomorphism
    since \(X\) as an object of \(\QSet_f\) is just \(i(X)\) for 
    \(X\) an object of \(\ScottQSet\), and hence adjunction 
    applies:
    \begin{align*}
        \QSet_f(\blank, X) 
        &= 
        \QSet_f(\blank, i(X)) 
        \\&\iso 
        \ScottQSet(\Singletons(\blank), X) 
        \\&=
        \QSet_f(i\circ\Singletons(\blank), i(X))
        \\&=
        \QSet_f(\Singletons(\blank), X)
        \intertext{
            Now, provided with the isomorphism but no guarantee 
            that \(X\) is complete, we have to find sufficient 
            reason for \(X\) to be complete.
        }
        \YonedaEmbedding_f(X) 
        &\iso
        \YonedaEmbedding_f(X)\circ i\circ\Singletons
        \\&=
        \ScottQSet_f(\Singletons(\blank), X)
        \\&\iso
        \ScottQSet_r(\Graph\circ\Singletons(\blank), \Graph(X))
        \\&=
        \ScottQSet_r(\Graph\circ\Singletons(\blank), X)
        \\&\iso
        \ScottQSet_r(\Graph(\blank), X)
        \\&=
        \YonedaEmbedding_r(X)\restriction\Graph
    \end{align*}        
\end{proof}

\begin{prop}
    The true force of the theorem above is that \(i\) in
    \(\Singletons\dashv i\) not only preserves limits but 
    indeed it creates them.
\end{prop}
\begin{proof}
    Take a diagram \(D\) on \(\ScottQSet\),
    \begin{align*}
        \YonedaEmbedding_f(\lim i\circ D) 
        &\iso
        \lim\YonedaEmbedding_f(i\circ D)
        \\&\iso
        \lim \YonedaEmbedding(i\circ D)\circ i\circ\Singletons
        \\&\iso
        \YonedaEmbedding(\lim i\circ D)\circ i\circ\Singletons
    \end{align*}
    This shows that the external (ie. in \(\QSet_f\)) the limit
    of complete \(\quantale Q\)-sets is itself complete. 
    
    Since any cone for \(D\) -- say, with vertex \(X\) -- would 
    \emph{be} a cone for \(i\circ D\) (as \(i\) is simply an 
    inclusion, they literally are) -- by universality there is
    exactly one morphism \(\lim i\circ D\to X\) making 
    everything commute.

    Since the external limit is now known to be complete, by 
    necessity, it also \emph{is} the limit in \(\ScottQSet\) as
    well. Hence, \(i\) both preserves all limits from its 
    domain, and reflects all limits in its codomain. 
\end{proof}
	\printbibliography{}
\end{document}